\documentclass[11pt]{article}

\usepackage{amscd,amsmath, amssymb, fancyhdr, epsfig,url,color}

\usepackage{soul}
\setstcolor{red}
\sethlcolor{red}

\numberwithin{equation}{section}

\newcommand{\version}{v. 1.11, March 31, 2022}

\def\eqref#1{(\ref{#1})}
\newcommand{\goth}{\mathfrak}

\newcommand{\Z}{{\mathbb Z}}
\newcommand{\C}{{\mathbb C}}

\newcommand{\R}{{\mathbb R}}
\newcommand{\Q}{{\mathbb Q}}

\def\1{\sqrt{-1}\:}

\newcommand{\cntrct}                
{\hspace{2pt}\raisebox{1pt}{\text{$\lrcorner$}}\hspace{2pt}}

\makeatletter
\def\x@arrow{\DOTSB\Relbar}
\def\xlongequalsignfill@{\arrowfill@\x@arrow\Relbar\x@arrow}
\newcommand{\xlongequal}[2][]{%
        \ext@arrow 0099\xlongequalsignfill@{#1}{#2}}
\def\xlongrightarrowfill@{\arrowfill@\relbar\relbar\longrightarrow}
\newcommand{\xlongrightarrow}[2][]{%
        \ext@arrow 0099\xlongrightarrowfill@{#1}{#2}}
\makeatother

\newcommand{\Lie}{\operatorname{Lie}}

\newcommand{\Diff}{\operatorname{Diff}}

\newcommand{\rk}{\operatorname{rk}}

\renewcommand{\Im}{\operatorname{Im}}

\newcounter{Mycounter}[section]
\newcounter{lemma}[section]
\renewcommand{\thelemma}{{Lemma \thesection.\arabic{lemma}}}
\newcommand{\lemma}{%
    \setcounter{lemma}{\value{Mycounter}}
    \refstepcounter{lemma}
    \stepcounter{Mycounter}
    {\noindent \bf \thelemma:\ }}

\newcounter{claim}[section]
\renewcommand{\theclaim}{{Claim \thesection.\arabic{claim}}}
\newcommand{\claim}{%
    \setcounter{claim}{\value{Mycounter}}
    \refstepcounter{claim}
    \stepcounter{Mycounter}
    {\noindent \bf \theclaim:\ }}

\newcounter{sublemma}[section]
\newcounter{corollary}[section]
\newcounter{theorem}[section]
\renewcommand{\thetheorem}{{Theorem \thesection.\arabic{theorem}}}
\newcommand{\theorem}{%
    \setcounter{theorem}{\value{Mycounter}}
    \refstepcounter{theorem}
    \stepcounter{Mycounter}
    {\noindent \bf \thetheorem:\ }}

\newcounter{conjecture}[section]
\newcounter{proposition}[section]
\renewcommand{\theproposition}
      {{Proposition \thesection.\arabic{proposition}}}
\newcommand{\proposition}{%
    \setcounter{proposition}{\value{Mycounter}}
    \refstepcounter{proposition}
    \stepcounter{Mycounter}
    {\noindent \bf \theproposition:\ }}

\newcounter{definition}[section]
\renewcommand{\thedefinition}
      {{Definition~\thesection.\arabic{definition}}}
\newcommand{\definition}{%
    \setcounter{definition}{\value{Mycounter}}
    \refstepcounter{definition}
    \stepcounter{Mycounter}
    {\noindent \bf \thedefinition:\ }}

\newcounter{example}[section]
\renewcommand{\theexample}{{Example \thesection.\arabic{example}}}
\newcommand{\example}{%
    \setcounter{example}{\value{Mycounter}}
    \refstepcounter{example}
    \stepcounter{Mycounter}
    {\noindent \bf \theexample:\ }}

\newcounter{remark}[section]
\renewcommand{\theremark}{{Remark \thesection.\arabic{remark}}}
\newcommand{\remark}{%
    \setcounter{remark}{\value{Mycounter}}
    \refstepcounter{remark}
    \stepcounter{Mycounter}
    {\noindent \bf \theremark:\ }}

\newcounter{problem}[section]
\newcounter{question}[section]
\makeatletter

\@addtoreset{equation}{section} \@addtoreset{footnote}{section} \makeatother

\def\proof{\hbox{\noindent {\bf Proof:}}}
\def\blacksquare{\hbox{\vrule width 5pt height 5pt depth 0pt}}
\def\endproof{\blacksquare}

\newcommand{\T}{{\mathbb{T}}}
\newcommand{\N}{{\mathbb{N}}}

\newcommand{\bGamma}{{\overline{\Gamma}}}

\newcommand{\bfU}{{\bf U}}

\newcommand{\cB}{{\mathcal{B}}}

\newcommand{\cL}{{\mathcal{L}}}

\newcommand{\cU}{{\mathcal{U}}}

\newcommand{\cX}{{\mathcal{X}}}

\newcommand{\hL}{{\hat{L}}}

\newcommand{\homega}{{\hat{\omega}}}

\newcommand{\Span}{{\rm Span}}

\newcommand{\tomega}{{\widetilde{\omega}}}

\newcommand{\Gu}{G}
\newcommand{\Hy}{H}
\newcommand{\Ky}{K}
\newcommand{\GLambda}{G_{\Lambda}}

\newcommand{\Subot}{{\mathcal{S}_{u^\bot}}}

\newcommand{\tSlin}{\cL}
\newcommand{\tSlinl}{\cL_l}
\newcommand{\tSlinlp}{\cL_{l'}}

\newcommand{\cLu}{\cL_u}

\newcommand{\TeichsK}{\operatorname{\sf T}}

\newcommand{\Sympl}{\operatorname{S}}
\newcommand{\SymplK}{\operatorname{\mathcal{S}_K}}
\newcommand{\SymplKone}{\operatorname{\mathcal{S}}}

\newcommand{\Pers}{\operatorname{\sf Per}}

\newcommand{\bigzero}{\mbox{\normalfont\Huge 0}}
\newcommand{\rvline}{\hspace*{-\arraycolsep}\vline\hspace*{-\arraycolsep}}

\begin{document}

\begin{center}
{\LARGE\bf
Rigidity of Lagrangian embeddings into symplectic tori and K3 surfaces\\
\ \\
}

Michael Entov,
Misha Verbitsky\footnote{
  Partially supported by  the  Russian Academic Excellence Project '5-100',
FAPERJ E-26/202.912/2018 and CNPq - Process 313608/2017-2.

 }

\end{center}

\begin{abstract}
\small
A K\"ahler-type form is a symplectic form compatible with
an integrable complex structure. Let $M$
 be either
 a
torus or a K3-surface equipped with a K\"ahler-type
form. We show that the homology class of any Maslov-zero
Lagrangian torus in $M$ has to be non-zero and
primitive. This extends previous results of Abouzaid-Smith
(for tori) and Sheridan-Smith (for K3-surfaces) who proved
it for particular K\"ahler-type forms on $M$. In the K3
case our proof uses dynamical properties of the action of
the diffeomorphism group of $M$ on the space of the
K\"ahler-type forms. These properties are obtained using
Shah's arithmetic version of Ratner's orbit closure theorem.
\end{abstract}


{\small
\tableofcontents
}


\section{Main results}
\label{_Intro_Main_Results_Subsection_}


Recall that
{\bf K\"ahler structure} on a smooth manifold $M$ is a pair $(\omega,J)$, where $\omega$ is a symplectic form on $M$ and $J$ is
an {\it integrable} complex structure on $M$ {\bf compatible} with
$\omega$, that is, satisfying $\omega(Jv, Jw)=\omega(v,w)$ and
$\omega(v,Jv) > 0$ for any non-zero tangent vectors $v,w$.
(We view complex structures as tensors on $M$ -- that is, as integrable almost complex structures; in particular, we only consider complex structures with the prefixed underlying smooth structure on $M$.) A symplectic form, or a complex structure, on $M$ is said to be of {\bf K\"ahler type}, if it appears in {\it some} K\"ahler structure.

In this paper we will be mostly concerned with $M$ which is either and even-dimensional torus $\T^{2n}=\R^{2n}/\Z^{2n}$ or a smooth manifold (of real dimension $4$) underlying a complex K3 surface.

Let us note the following:

\bigskip
\noindent
- The K\"ahler-type symplectic forms on $\T^{2n}$ are exactly the ones that can be mapped
by an orientation-preserving diffeomorphism of $\T^{2n}$ to a linear symplectic form
-- see e.g. \cite[Prop. 6.1]{_EV-JTA_}. (A linear symplectic form on $\T^{2n}$ is a form whose lift to $\R^{2n}$
has constant coefficients). We
will always use the orientation on $\T^{2n}$ induced by
the standard orientation on $\R^{2n}=\C^n$ and consider only the K\"ahler-type symplectic forms on $\T^{2n}$ inducing this orientation.

\bigskip
\noindent
- All smooth manifolds $M$ underlying a complex K3 surface
are diffeomorphic. They all are compact and connected; any
complex structure on such an $M$ is of K\"ahler-type and
its first Chern class is zero. In fact, all complex
structures
 and K\"ahler-type symplectic forms on such an $M$ appear as parts of
 hyperk\"ahler
structures
 -- see
\ref{_Kahler-type-on-K3-is-hyperkahler-type_Proposition_}
below. All complex structures on such an $M$ define the
same orientation on $M$ so that $b_+ (M)=3$ and
${b_-} (M)=19$ -- see e.g. \cite{_Geom-K3-Asterisque_} (reversing the
orientation would yield $b_+ =19$ and $b_-=3$).
Further on, we will always equip a K3 surface with this
standard orientation. Alternatively, to fix the orientation,
one can use the fact that all complex K3 surfaces are deformation equivalent and, in
particular, oriented diffeomorphic -- see e.g. \cite{_Geom-K3-Asterisque_}.

\hfill

This is the main application of the
 techniques
 developed in
this paper.


\hfill


\theorem\label{_Main_Theorem_}

Assume $M$, $\dim_\R M = 2n$, is either an even-dimensional torus or a smooth manifold underlying a K3 surface. Let
$\omega$ be a K\"ahler-type symplectic form on $M$.

Then for any Maslov-zero Lagrangian torus $L\subset (M, \omega)$ the homology class $[L] \in H_n (M;\Z)$ is non-zero
and primitive.


\hfill


For the proof see Section~\ref{_pfs-main-thm_Section_}.

In the torus case the result of \ref{_Main_Theorem_} was previously proved by Abouzaid and Smith \cite[Cor. 1.6 and Cor. 9.2]{_Abouzaid-Smith_}
for {\it one particular} K\"ahler-type symplectic form on
$M=\T^{2n}$ -- namely, for the standard Darboux form. Their proof of $[L]\neq 0$ relies on a previous deep
result of Fukaya and works similarly for any linear (or, equivalently, K\"ahler-type) symplectic form on $\T^{2n}$ -- we recall their argument in the proof of \ref{_Main_Theorem_}.

In the K3 case the result of \ref{_Main_Theorem_} was previously proved by Sheridan and Smith \cite[Thm. 1.3]{_Sheridan-Smith_} for {\it some} K\"ahler-type symplectic forms on $M$. Their proof
uses deep methods of homological mirror symmetry.

Here we extend the results of Abouzaid-Smith (in the torus case) and Sheridan-Smith (in the K3 case) to {\it any} K\"ahler-type symplectic form on $M$. Our result in the K3 case answers a question of Seidel (see \cite[Sect. 1.1, Quest. 5]{_Sheridan-Smith_}).

Our proof of \ref{_Main_Theorem_} reduces the case of an arbitrary K\"ahler-type symplectic form on $M$ to particular cases where symplectic rigidity results can be applied. Namely, if $(M,\omega)$ admits a Maslov-zero Lagrangian torus $L$ whose homology class is zero or non-primitive, then so does $(M,\omega')$ for any symplectic form $\omega'$ sufficiently close to $\omega$ and satisfying $[\omega'|_L]=0$ (this is an easy application of Moser's method). We show that such $\omega'$ can be chosen to lie in the $\Diff^+ (M)$-orbit of a K\"ahler-type symplectic form for which the absence of such Lagrangian tori is proved by symplectic rigidity methods.
 (Here and further on $\Diff^+ (M)$ stands for the group of orientation-preserving diffeomorphisms of $M$.)
 Since the set of the symplectic forms on $M$ admitting such a Lagrangian torus is $\Diff^+ (M)$-invariant, this yields a contradiction that proves \ref{_Main_Theorem_}.

In the torus case the application of symplectic rigidity results follows closely the proof of Abouzaid-Smith in \cite{_Abouzaid-Smith_} and in the K3 case we use the work of Sheridan and Smith \cite{_Sheridan-Smith_}.

The dynamical result on the $\Diff^+ (M)$-orbits needed to prove $[L]\neq 0$ in the K3 case can be deduced relatively easily from \cite{_EV-JTA_}, and the main effort in this paper is to strengthen it in order to prove the
 {\it primitivity}
 of $[L]$.


\hfill


\remark

The existence of Maslov-zero Lagrangian tori is well-known for certain K\"ahler-type symplectic forms on tori and K3 surfaces. For instance,

\begin{itemize}

\item{} For the standard Darboux form $dp\wedge dq$ on a torus, the meridian Lagrangian torus $\{ p = \textrm{const}\}$ is Maslov-zero.

\item{} K3 surfaces with certain K\"ahler forms admit special Lagrangian tori \cite{_Harvey-Lawson_} (also see \cite{_SYZ_}). These special Lagrangian tori are Maslov-zero.

\end{itemize}

At the same time there exist K\"ahler-type symplectic forms on tori and K3 surfaces that do not admit any Lagrangian submanifolds with a non-zero homology class (for instance, because of obvious homological obstructions), and, in particular, by \ref{_Main_Theorem_}, no Maslov-zero Lagrangian tori.


\hfill


\remark

Using the Kodaira-Spencer stability theorem
\cite{_Kod-Spen-AnnMath-1960_} and Torelli theorems for even-dimensional
tori (see e.g. \cite[Ch. I,
  Thm. 14.2]{_Barth-Hulek-Peters-vdVen_}) and K3 surfaces
(see \cite{_Burns-Rapoport_},
cf. \cite[p.96]{_Geom-K3-Asterisque_}) one can show that
the K\"ahler-type symplectic forms on these manifolds form
an open subset of the space of all symplectic forms (with
respect to the $C^\infty$-topology).

It is an open question whether {\it all} symplectic forms on
$\T^{2n}$ and K3 surfaces
 are of K\"ahler type. In
\cite{_Donaldson:ellipt_} Donaldson outlined how one can
try to prove that the answer to the question is {\it
  positive}.
On other hand, in view of \ref{_Main_Theorem_}, a possible strategy of proving that the answer to the question is {\it negative} would be to construct a symplectic form on $\T^{2n}$, or on a K3 surface, that admits
a Maslov-zero
Lagrangian torus whose integral homology class is zero or non-primitive.


\hfill


Let us now describe the dynamical results at the heart of the proof of \ref{_Main_Theorem_} in more detail.

First, let us consider the case of $\T^{2n}$.

Consider the linear symplectic forms of volume $1$ on $\T^{2n}$ and denote by
$\tSlin$ the space of their lifts to $\R^{2n}$. The
topology on $\tSlin$ is induced by the standard topology
on the space of bilinear forms on $\R^{2n}$.

Let $l\subset \R^{2n}$ be a vector subspace of (real) dimension $n$ which is spanned over $\R$ by vectors in $\Z^{2n}\subset \R^{2n}$ -- such a vector subspace is called {\bf rational}.
Define $\tSlinl\subset \tSlin$ as
\[
\tSlinl := \{\ \tomega\in\tSlin\ |\ \tomega|_l\equiv 0\ \}.
\]
Define a group $G\subset SL (2n,\R)$ by
\[
G: = \{\ g\in SL (2n,\R)\ |\ g|_l= Id\ \}.
\]
The group $G$ acts on $\tSlinl$.
Let
\[
\GLambda := G\cap SL (2n,\Z).
\]

We say that $n$-dimensional rational subspaces $l,l'\subset\R^{2n}$ are {\bf complementary} if
$(l\cap \Z^{2n})\oplus (l'\cap\Z^{2n}) = \Z^{2n}$ (in particular, this implies that $l\oplus l' = \R^{2n}$).

Given complementary $n$-dimensional rational subspaces $l,l'\subset\R^{2n}$, we say that a symplectic form $\tomega\in \tSlin$ is {\bf $(l,l')$-Lagrangian split} if $l$ and $l'$ are Lagrangian with respect to $\tomega$, that is, $\tomega\in \tSlinl\cap\tSlinlp$.


\hfill


\proposition\label{_GLambda-orbit-torus-case_Proposition_}

The union of the $G_\Lambda$-orbits of the $(l,l')$-Lagrangian split forms
is dense in $\tSlinl$.


\hfill


For the proof see Section~\ref{_linear-sympl-forms-on-tori_Section_}.

Now let $M$ be either $\T^4$ or a smooth manifold underlying a K3 surface.

Denote by $\SymplKone (M)$ the space of K\"ahler-type symplectic forms $\omega$ on $M$ of total volume $1$, meaning that $\int_M \omega^n = 1$, $2n = \dim_R M$.
We equip $\SymplKone (M)$ with the $C^\infty$-topology.

Let $(\cdot,\cdot): H^2 (M;\R)\times H^2 (M;\R)\to\R$, $(x,y) := \langle x\cup y, [M]\rangle$, be the intersection product.

\bigskip
\noindent
{\bf Further on, the positive/unit/isotropic vectors in $H^2 (M;\R)$, as well the orthogonal complement $x^\bot$ of $x\in H^2 (M;\R)$,
are all taken with respect to the bilinear symmetric form $(\cdot,\cdot)$.}


\hfill


\definition\label{_isotr-orthoirrational_Definition_}

We say that a positive $y\in H^2 (M;\R)$ is
{\bf orthoisotropically irrational} if

\smallskip
\noindent
(A) There exists an isotropic $u\in H^2 (M;\Z)$ such that $y\in u^\bot$.

\smallskip
\noindent
(B) For any $u\in H^2 (M;\Z)$ such that $y\in u^\bot$
one has
$\Span_\R \{ u,y\}\cap H^2 (M;\Z) = \Span_\Z \{ u\}$. (In particular, $y$ is a not a multiple of an integral cohomology class).

\hfill

\remark

The set ${\goth S}$ of orthoisotropically irrational classes
be\-longs to the union $\bigcup_{u\in H^2 (M;\Z), (u,u)=0}u^\bot$
of countably many orthogonal complements
$u^\bot$ for isotropic $u\in H^2 (M;\Z)$.
It can be obtained from the union above
by re\-mo\-ving countably many co\-dimen\-sion-2 sub\-spaces:
\[
{\goth S}= \bigcup_{u\in  H^2 (M;\Z), (u,u)=0}\left (u^\bot
\setminus\ \bigcup_{v\in H^2 (M;\Z)\cap u^\bot, \ \dim
  \Span \{ u, v\} =2} \Span \{ u, v\} \right).
\]


\hfill


Let $P\subset H^2 (M;\R)$ be the set of unit vectors:
\[
P:= \{\ y\in H^2 (M;\R)\ |\ (y,y)=1\ \}.
\]
For $u\in H^2 (M;\Z)$ define a subspace $\Subot (M)\subset \SymplKone (M)$ by
\[
\Subot (M):= \{\ \omega\in\ \SymplKone (M)\ |\ [\omega]\in u^\bot\cap P\ \}.
\]

The group $\Diff^+ (M)$ acts naturally on $\SymplKone (M)$.


\hfill


\theorem\label{_dense-orbit-of-sympl-form-K3-case_Theorem_}

Assume that $\omega_0\in \SymplKone (M)$ so that the cohomology class $[\omega_0]\in H^2 (M;\R)$ is orthoisotropically irrational.

Then for any isotropic $u\in H^2 (M;\Z)$
the intersection of the $\Diff^+ (M)$-orbit of $\omega_0$ with $\Subot (M)$
is dense in $\Subot (M)$. In particular, it is dense in $\SymplKone (M)$ (for $u=0$).


\hfill


For the proof of
\ref{_dense-orbit-of-sympl-form-K3-case_Theorem_} see
Section~\ref{_Sympl-Teichm-space-mapp-class-gr_Subsection_}. It
follows ideas similar to \cite{_EV-JTA_} and going back to
\cite{_V-Duke_} and is based on Shah's arithmetic version
\cite{_Shah:uniformly_} of the famous Ratner's orbit closure theorem
\cite{_Ratner_Duke_1991_}.

Let us compare \ref{_dense-orbit-of-sympl-form-K3-case_Theorem_} to the results in \cite{_EV-JTA_}.

In the case $M=\T^4$ it was proved in \cite{_EV-JTA_} that for any K\"ahler-type symplectic form $\omega_0\in\SymplKone (M)$,
such that $[\omega_0]\in H^2 (M;\R)$ is not a real multiple of a rational cohomology class, the $\Diff^+ (M)$-orbit of $\omega_0$
is dense in $\SymplKone (M)$.

In the K3 case it was proved in \cite{_EV-JTA_} that for
any K\"ahler-type\footnote{In \cite{_EV-JTA_} we used the term ``hyperk\"ahler type''
for a symplectic form which can be obtained
as a K\"ahler form of some hyperk\"ahler structure.
Using Moser's lemma and Calabi-Yau theorem, is easy to see
that on a K3 surface
any symplectic form of K\"ahler-type is, in fact, of hyperk\"ahler type (Section~\ref{_K3-surfaces_Section_}).} symplectic form
$\omega_0\in\SymplKone (M)$,
such that $[\omega_0]\in H^2 (M;\R)$ is not a real
multiple of a rational cohomology class, the $\Diff^+
(M)$-orbit of $\omega_0$
is dense in the connected component of $\omega_0$ in the
space of K\"ahler-type symplectic forms lying
in $\SymplKone (M)$.

Clearly, if $[\omega_0]$ is orthoisotropically irrational, then it is not a real multiple of a rational cohomology class.
Thus, \ref{_dense-orbit-of-sympl-form-K3-case_Theorem_} strengthens the results in \cite{_EV-JTA_}: it shows that for any $\omega_0\in \SymplKone (M)$ such that $[\omega_0]$ is orthoisotropically irrational the $\Diff^+ (M)$-orbit of $\omega_0$ is not only dense $\SymplKone (M)$
(which follows from \cite{_EV-JTA_}) but also that its intersection with $\Subot (M)$ is dense in $\Subot (M)$ for any isotropic $u\in H^2 (M;\Z)$.

Let us now state a proposition which is used in the proof of the results above and may be of independent interest.


\hfill


\proposition
\label{_DiffH-acts-transitively-on-cohomologous-Kahler-type-forms_Proposition_}

Let $M^{2n}$ be either $\T^{2n}$ or a smooth manifold underlying a complex K3 surface.

Then any two K\"ahler-type symplectic forms on $M$ (compatible with the orientation of $M$) can be mapped
into each other by a diffeomorphism of $M$ acting trivially on homology.


\hfill


For the proof see Section~\ref{_Sympl-Teichm-space-mapp-class-gr_Subsection_}.


\hfill


\section{Linear symplectic forms on tori -- proof of \ref{_GLambda-orbit-torus-case_Proposition_}}
\label{_linear-sympl-forms-on-tori_Section_}

Let us prove \ref{_GLambda-orbit-torus-case_Proposition_}.

Since $l$ and $l'$ are complementary, one can choose a basis $\cB$ of $\R^{2n}$ formed by vectors in $\Z^{2n}$
so that the first $n$ basic vectors lie in $l$ and the last $n$ ones lie in $l'$.

Let $M_n (\R)$ (respectively $M_n (\Z)$) denote the spaces of $n\times n$-matrices with real (respectively integral) coefficients.

With respect to the basis $\cB$:

\smallskip
\noindent
- The matrices of the elements of $G$  are exactly the matrices of the form
\begin{equation}
\label{_torus-case-G-matrices_Equation_}
\begin{pmatrix}
I_n & B\\
0 & A
\end{pmatrix},
\end{equation}
where
$A\in SL (n,\R)$, $B\in M_n (\R)$.

\smallskip
\noindent
- The matrices of the elements of $G_\Lambda$  are exactly the matrices of the form \eqref{_torus-case-G-matrices_Equation_} with integral coefficients.

\smallskip
\noindent
- The matrices of the forms in $\tSlinl$ are exactly the matrices of the form
\[
\begin{pmatrix}
0 & -C^t\\
C & D
\end{pmatrix},
\]
where $C\in SL (n,\R)$ and $D$ is skew-symmetric.

\smallskip
\noindent
- The matrices of the
$(l,l')$-Lagrangian split forms are exactly the matrices of the form
\[
\begin{pmatrix}
0 & -C^t\\
C & 0
\end{pmatrix},
\]
where $C\in SL (n,\R)$.

We will identify the elements of $G$ and the linear symplectic forms with the corresponding matrices.

Define $\cX\subset SL (n,\R)$ by
\[
\cX := \left\{\ C\in SL (n,\R)\ |\ \textrm{the entries of}\ C^{-1}\ \textrm{are linearly independent over}\ \Q\ \right\}.
\]
The set $\cX$ satisfies the following properties:

\smallskip
\noindent
(I) $\cX$ is dense in $SL (n,\R)$.

\smallskip
\noindent
(II) For each $C\in\cX$ the projection of the set $\R C^{-1}$ to the torus $T:=M_n (\R)/M_n (\Z)$ is dense in $T$.

\smallskip
Property (II) implies the following claim:

\smallskip
\noindent
(III) For each $C\in\cX$ the set $\{ CB - (CB)^t, B\in M_n (\Z)\}$ is dense in the space of skew-symmetric $n\times n$-matrices.

\smallskip
Indeed, by (II), for each skew-symmetric $n\times n$-matrix $D$ the matrix $C^{-1} D$ can be approximated by matrices of the form $tC^{-1}+B$, $t\in\R$, $B\in M_n (\Z)$. Accordingly, $D$ can be approximated by matrices $tId + CB$, $t\in\R$, $B\in M_n (\Z)$. Since any matrix can be
uniquely represented as a sum of a symmetric and a skew-symmetric matrices, this means that $D$ can be approximated by the skew-symmetric components of the matrices $CB$, $B\in M_n (\Z)$ -- that is, by the matrices $CB - (CB)^t$, $B\in M_n (\Z)$, which proves (III).

Define $X\subset \tSlinl$ as the set of the $(l,l')$-Lagrangian split forms represented by the matrices
\[
\begin{pmatrix}
0 & -C^t\\
C & 0
\end{pmatrix},
C\in\cX.
\]

In order to prove the proposition it suffices to show that
the union of the $G_\Lambda$-orbits of the forms in $X$ is dense in $\tSlinl$ -- let us prove this claim.

The action of the element
\[
\begin{pmatrix}
I_n & B\\
0 & I_n
\end{pmatrix}\in G
\]
on a form
\[
\begin{pmatrix}
0 & -C^t\\
C & 0
\end{pmatrix}\in X
\]
yields the form
\[
\begin{pmatrix}
I_n & B\\
0 & I_n
\end{pmatrix}^t
\begin{pmatrix}
0 & -C^t\\
C & 0
\end{pmatrix}
\begin{pmatrix}
I_n & B\\
0 & I_n
\end{pmatrix}
=
\]
\[
=
\begin{pmatrix}
I_n & 0\\
B^t & I_n
\end{pmatrix}
\begin{pmatrix}
0 & -C^t\\
C & 0
\end{pmatrix}
\begin{pmatrix}
I_n & B\\
0 & I_n
\end{pmatrix}
=
\begin{pmatrix}
0 & -C^t\\
C & C B - B^t C^t
\end{pmatrix}.
\]
By (I) and (III) above, the set of the matrices
\[
\begin{pmatrix}
0 & -C^t\\
C & C B - B^t C^t
\end{pmatrix},\ C\in\cX,\ B\in M_n (\Z),
\]
is dense in the set of the matrices
\[
\begin{pmatrix}
0 & -C^t\\
C & D
\end{pmatrix}, C\in SL (n,\R),\ D\ \textrm{is skew-symmetric}.
\]
In other words,
 any
 form in $\tSlinl$
can be approximated by the image of
a form
\[
\begin{pmatrix}
0 & -C^t\\
C & 0
\end{pmatrix}\in X
\]
under the action of
\[
\begin{pmatrix}
I_n & B\\
0 & I_n
\end{pmatrix}\in G_\Lambda
\]
for appropriate $C\in\cX$ and $B\in M_n (\Z)$.
This proves the claim and the proposition.
\endproof


\hfill


\section{Lattices}
\label{_Lattices-Ratner_Section_}

Let us start with recalling a few generalities on lattices and Ratner's orbit closure theorem.


\hfill


\subsection{Quadratic lattices}
\label{_Quadratic lattices_Subsection_}

A {\bf quadratic vector space} is a finite-dimensional real vector space equip\-ped with an $\R$-valued symmetric product $(\cdot,\cdot)$, bilinear over $\R$.

A {\bf quadratic (or integral) lattice} is a free finite-rank abelian group equipped with a $\Z$-valued symmetric product $(\cdot,\cdot)$, bilinear over $\Z$.

Given a quadratic vector space $V$, a {\bf quadratic lattice in $V$} is a discrete subgroup $\Lambda\subset V$ which together with the restriction of the product from $V$ forms a quadratic lattice (that is, the restriction of the product to $\Lambda$ is $\Z$-valued).

Given a quadratic lattice $\Lambda$, the $\Z$-valued product on $\Lambda$ extends to an $\R$-valued bilinear symmetric product on the real vector space $V:=\Lambda\otimes_\R \R$, that will be also denoted by $(\cdot,\cdot)$, and $\Lambda$ becomes a quadratic lattice in $V$.

\bigskip
{\bf
Further on in Section~\ref{_Quadratic lattices_Subsection_}, $V$ will stand for $\Lambda\otimes_\R \R$ and we will view the elements of $\Lambda$ as vectors in $V$.
}
\bigskip

A quadratic lattice $\Lambda$ is called {\bf unimodular} if the matrix of $(\cdot,\cdot)$ with respect to a basis of $\Lambda$ (over $\Z$) has determinant $\pm 1$. It is called {\bf even} if the number $(x,x)$ is even for all $x\in\Lambda$.

The {\bf signature of $\Lambda$} is defined as the signature of the symmetric bilinear form $(\cdot,\cdot)$ on $\Lambda\otimes_\R \R$.


\hfill


\example
\label{_integral hyperbolic plane_Example_}

The {\bf quadratic (or integral) hyperbolic plane} $\bfU$ is a quadratic lattice of signature $(1,1)$: $\bfU := \Span_\Z \{ x,y\}$, $(x,x)=(y,y)=0$, $(x,y)=(y,x)=1$. This quadratic lattice is even and unimodular.


\hfill


The notions of a quadratic sublattice and the direct sum of quadratic lattices are defined in an obvious way.

A vector $x \in \Lambda$ is called {\bf primitive} if it is not an integral multiple of another element in $\Lambda$.

A vector $x \in V$ is called {\bf isotropic} if $(x,x)=0$, {\bf positive} if $(x,x)>0$, and {\bf unit} if $(x,x)=1$.

The notation $^\bot$ will be used for orthogonal complements in $V$ with respect to $(\cdot,\cdot)$.

Denote by $x^{\bot,1}$ the set of unit vectors in $x^\bot$:
\[
x^{\bot,1} := \{\ y\in x^\bot\ |\ (y,y)=1\ \}.
\]


\hfill


\proposition\label{_unimodular-lattices-divisor-of-any-integral-vector-is-one_Proposition_}

Assume $\Lambda$ is a unimodular quadratic lattice and $u\in\Lambda$ is a non-zero vector.

Then the following claims hold:

\smallskip
\noindent
(A) If $u$ is primitive, then $(u,\Lambda)=\Z$.

\smallskip
\noindent
(B) Assume $\Lambda$ is even and of signature $(p,q)$, and $u$ is isotropic. Then

\smallskip
\noindent
(B1) $u^\bot \cap \Lambda = \Span_\Z \{ u\} \oplus \Lambda'$, where $\Lambda'$ is a unimodular integral lattice in $u^\bot$ of signature $(p-1,q-1)$;

\smallskip
\noindent
(B2) For any $x\in u^\bot \cap \Lambda$, such that $x\notin\Span_\Z \{ u\}$, the rank of the group $(\Span_\R \{ u,x\})^\bot \cap \Lambda$ is $\rk \Lambda - 2$.


\hfill


\noindent
{\bf Proof of \ref{_unimodular-lattices-divisor-of-any-integral-vector-is-one_Proposition_}.}

Part (A) is an elementary algebraic exercise based on the fact that the ideal in $\Z$ generated by a finite collection of integers is a principal ideal generated by their g.c.d.

Let us prove part (B). Without loss of generality, we may assume that $u$ is primitive.

By part (A), there exists a vector $s\in\Lambda$ such that $(u,s)=1$. Then, since $\Lambda$ is even, $(s,s)$ is even and $z:= s - (s,s)u/2$ is an isotropic element of $\Lambda$ that spans together with $u$ an integral hyperbolic plane $\bfU$. Since $\bfU$ is a unimodular sublattice of $\Lambda$, we can write $\Lambda$ as $\Lambda=\bfU\oplus (\bfU^\bot \cap \Lambda)$ (see e.g. \cite[Ch.1, Lem.3.1]{_Milnor-Huse_}). The signature of $\bfU^\bot$ is $(p-1,q-1)$. Setting $\Lambda' := \bfU^\bot \cap \Lambda$ finishes the proof of (B1).

Claim (B2) follows immediately from the following general claim: if $\Lambda_1$ is a sublattice of $\Lambda$, then $\Lambda_1^\bot := (\Lambda_1 \otimes_\Z \Q)^\bot\cap\Lambda$ has rank $\rk \Lambda - \rk \Lambda_1$. (The claim follows from the observation that any vector in $(\Lambda_1 \otimes_\Z \Q)^\bot$ has an integral multiple that lies in $\Lambda_1^\bot$).
\endproof


\hfill


Denote by $SO (V)$ the group of linear automorphisms of $V$, with determinant $1$, preserving the bilinear form. Let
$SO^+ (V)$ be the identity component of $SO (V)$, and let $SO^+ (\Lambda)$ be the subgroup of $SO^+ (V)$ formed by the elements of $SO^+ (V)$ preserving $\Lambda$.


\hfill


\definition
\label{_isotropically-transitive-lattice_Definition_}

We say that a quadratic lattice $\Lambda$ is {\bf isotropically transitive} if $SO^+ (\Lambda)$ acts transitively on the set of primitive isotropic vectors in $\Lambda$.


\hfill


\example
\label{_U-U-isotropically transitive_Example_}

Assume that $\Lambda = \bfU\oplus\bfU \oplus\Lambda'$, where $\Lambda'$ is an even unimodular quadratic lattice of signature $(p,q)$, $p,q\in\N$.

Then $\Lambda$ is isotropically transitive by
\cite[Prop.3.3(i)]{_Gritsenko-et-al_} (part (A) of
\ref{_unimodular-lattices-divisor-of-any-integral-vector-is-one_Proposition_}
allows to apply the latter result).

\hfill

\remark

We will apply this result to the set of isotropic cohomology
classes on a K3 surface (Section~\ref{_density-of-Diff-orbits-in-cohomology_Subsection_}).
We will get that, since  the cohomology
classes of Lagrangian tori are isotropic, the
mapping class group acts transitively on the lines spanned by these classes
in $H^2(M; \Q)$.


\hfill


\definition
\label{_isotropically-irrational-vector_Definition_}

Let $y\in V$ be a positive vector.

Given a non-zero isotropic $u\in\Lambda$, we say that $y$
is {\bf $u$-orthoirrational} if $y\in u^\bot$ and $y$ does
not lie in a two-dimensional plane spanned by $u$ and a
vector in $\Lambda\cap u^\bot$.

We say that $y$ is {\bf orthoisotropically irrational} if
if it lies in $u^\bot$ for some non-zero isotropic
$u\in\Lambda$ and is $u$-orthoirrational with respect to
any non-zero isotropic $u\in y^\bot\cap\Lambda$.
(This matches \ref{_isotr-orthoirrational_Definition_}).


\hfill


\subsection{Lattices in Lie groups and Shah's version of Ratner's orbit closure theorem}
\label{_Lattices-in-Lie-groups_Subsection_}

Let $G$ be a Lie group.

The group $G$ admits a left-invariant (respectively, right-invariant) measure, called {\bf the left (respectively, right) Haar measure}, which is defined uniquely, up to a constant factor -- this measure is defined by a left-invariant (respectively, right-invariant) differential volume form on $G$.

Assume $\Gamma\subset G$ is a discrete subgroup. The restriction of the right Haar measure on $G$ to
 the
 fundamental domain of the right action of $\Gamma$ on $G$ induces a measure $\mu_{G/\Gamma}$ on the space $G/\Gamma$ (this measure can be defined by a differential volume form on $G/\Gamma$).

The subgroup $\Gamma$ is called a {\bf Lie lattice}\footnote{This meaning of the term ``lattice" is not to be confused with the one in Section~\ref{_Quadratic lattices_Subsection_}.} (in $G$) if $\mu_{G/\Gamma} (G/\Gamma)< +\infty$.

Recall that $g\in G$ is called {\bf
unipotent}, if $g=e^h$ for an $ad$-nilpotent element $h$ of the Lie algebra
of $G$.
A {\bf unipotent
 one-parameter
 subgroup of $G$} is a subgroup of the form $\{ e^{th}\}_{t\in\R}$ for an $ad$-nilpotent element $h$ of the Lie algebra of $G$.
We say that $G$ is {\bf generated by unipotents} if
it is multiplicatively generated by unipotent elements.

A {\bf linear} (real) Lie group is a Lie subgroup of $SL (V)$ for a finite-dimensional real vector space $V$.

Let $V$ be a finite-dimensional real vector space, $m=\dim_\R V$, and let $\Lambda\subset V$ be a free discrete subgroup of $V$ of rank $m$. A basis of $V$ is called {\bf integral}, if it is a basis of $\Lambda$ over $\Z$.

Let $G\subset SL (V)$ be a linear Lie group.

Let
\[
\GLambda := \{\ g\in G\ |\ g(\Lambda)=\Lambda\ \}.
\]
A Lie lattice
$\Gamma$ in $G$ is called {\bf arithmetic} if $\Gamma \cap
\GLambda$ has a finite index in $\Gamma$ and $\GLambda$.

A {\bf $\Q$-character} on $G$ is a homomorphism $G\to
 \R_{>0}$ which is defined by algebraic equations with rational
 coefficients on the entries of the real $m\times
 m$-matrices representing the elements of $G$ with respect
 to an integral basis of $V$.

We say that $G$ is an {\bf algebraic group}
(respectively, an {\bf algebraic $\Q$-group})
if for an integral basis of $V$ the $m\times m$-matrices
of the elements of $G$ with respect to the basis
form a Lie subgroup of
$SL (m,\R)$ defined by algebraic equations with real
(respectively, rational)
coefficients on the entries of the
matrices.


\hfill


\proposition
\label{_SO-p-q-generated-by-unipotents_Proposition_}

Assume $p,q\in\Z_{>0}$ and $p+q>2$. Then the group $SO^+
(p,q)$ (the identity component of $SO(p,q)$)
is generated by its algebraic
 one-parameter
 unipotent subgroups.


\hfill


\proof
The group $SL(2,\R)$
contains a non-trivial algebraic
 one-parameter
 uni\-potent subgroup
\[
U:=\left\{
\begin{pmatrix}
1 & t \\ 0 & 1
\end{pmatrix},
t\in\R
\right\}.
\]
There is an isomorphism between $SL(2,\R)/\{\pm 1\}$ and $SO^+(1,2)$ defined by the following homomorphism $\mu: SL(2,\R)\to SO^+(1,2)$
with the kernel $\{\pm 1\}$: identify $\R^3$ with the space $V$ of real trace-free $2\times 2$-matrices equipped with the bilinear symmetric form $\langle A,B\rangle\to tr (AB)$, $A,B\in V$, of signature $(1,2)$; for each $C\in SL(2,\R)$ define $\mu(C)$ as the automorphism of $V$ given by $A\mapsto CAC^{-1}$, $A\in V$. The homomorphism $\mu$ is a polynomial map from $SL (2,\R)$ to $SO (1,2)$. Hence, $\mu (U)$ is a non-trivial algebraic
 one-parameter
 unipotent subgroup of $SO^+(1,2)$. This implies that for all $p,q\in\Z_{>0}$ and $p+q>2$ the group $SO^+ (p,q)$ also contains non-trivial algebraic
 one-parameter
 unipotent subgroups.

Any conjugate of an algebraic
 one-parameter
 unipotent subgroup of a real Lie group is again an algebraic
 one-parameter
 unipotent subgroup of the same Lie group. Therefore
the subgroups of $SO^+ (p,q)$ and of $SL(2,\R)$ generated by their algebraic
 one-parameter
 unipotent subgroups are normal and non-discrete.
Unless $p=q=2$, the group $SO^+ (p,q)$ is a simple Lie group and so is $SL(2,\R)$. Therefore these groups do not contain non-discrete proper normal subgroups \cite{_Ragozin_}. Consequently, they are generated by their algebraic
 one-parameter
 unipotent subgroups.

Let us consider the remaining case $p=q=2$. We have
$SO^+ (2,2) = \big(SL (2,\R)\times SL (2,\R)\big)/\{\pm 1\}$.
If $U$ is an arbitrary algebraic
 one-pa\-ra\-me\-ter
 uni\-po\-tent subgroup of $SL (2,\R)$, then
$Id\times U$ and $U\times Id$ are algebraic
 one-pa\-ra\-me\-ter
 unipotent subgroups of $SO^+ (2,2) = \big(SL (2,\R)\times SL (2,\R)\big)/\{\pm 1\}$.
Consequently, since $SL(2,\R)$ is generated by its algebraic
 one-pa\-ra\-me\-ter
 unipotent subgroups, so is $SO^+ (2,2)$.

This shows that $SO^+ (p,q)$ is generated by algebraic
 one-parameter
uni\-potent subgroups for all $p,q\in\Z_{>0}$ and $p+q>2$.
\endproof


\hfill


The following theorem belongs to A.Borel and Harish-Chandra.


\hfill


\claim (Borel-Harish-Chandra theorem, \cite[Thm. 9.4]{_Borel-HC1962_})
\label{_Borel-Harish-Chandra_Claim_}

Let $G\subset SL(V)$ be an algebraic $\Q$-group. Then
$\GLambda$ is a Lie lattice in $G$ (or, equivalently, $G$ admits an arithmetic lattice)
if and only if $G$ does not admit non-trivial $\Q$-characters.
\endproof


\hfill


\proposition
\label{_lattice-identity-component_Proposition_}

Let $G'\subset SL(V)$ be an algebraic $\Q$-group and $G$ its identity component (in the Lie group topology). Assume
$G$ does not admit non-trivial $\Q$-characters.

Then $G_\Lambda$ is a Lie lattice in $G$.


\hfill


\proof

Since $G$ does not admit non-trivial $\Q$-characters, $G'$ does not admit them either. Therefore, by Borel and Harish-Chandra theorem
(\ref{_Borel-Harish-Chandra_Claim_}), $G'_\Lambda$ is a Lie lattice in $G'$.

The identity component $G$ is a normal subgroup of $G'$ and its index is finite, since, by Whitney's theorem \cite{_Whitney_}, the real algebraic
subvariety of a Euclidean space has finitely many connected components. This implies that $G_\Lambda$ has a finite index in $G'_\Lambda$.
Since the indices $|G':G|$ and $|G'_\Lambda:G_\Lambda|$ are finite and $G'_\Lambda$ is a Lie lattice in $G'$, one easily gets that $G_\Lambda$ is a Lie lattice in $G$.
\endproof


\hfill


The following fundamental theorem belongs to M.Ratner.


\hfill


\claim (Ratner's orbit closure theorem, \cite{_Ratner_Duke_1991_})
\label{_Ratner_orbit_closure_Claim_}

Let $G$ be a
connected
Lie group,  $H\subset G$ its Lie subgroup
generated by unipotents and $\Gamma\subset G$ a Lie lattice.
Then for any $g\in G$ the closure $\overline{\Gamma g H}$
of the double class is obtained as
\[
\overline{\Gamma g H} = \Gamma g S
\]
for some closed Lie subgroup $S$, $H\subset S\subset G$.
In particular, if $H$ is a closed Lie subgroup, then the closure of
the orbit $\Gamma \cdot g H$ of $gH$ in
$G/H$ is $\Gamma (gSg^{-1}) \cdot gH$.
\endproof


\hfill


A combination of Ratner's orbit closure theorem (\ref{_Ratner_orbit_closure_Claim_}) with
a result of Shah \cite[Proposition 3.2]{_Shah:uniformly_}) yields a
more precise description of the group $S$ in the case where $G$ is
the identity component of
a linear algebraic $\Q$-group and $\Gamma\subset G$ is an arithmetic
Lie lattice.
We will state the result for $g=e$, since this is exactly what we
are going to use in our proof.


\hfill


\claim \label{_Ratner_arithmetic_Claim_}
(\cite[Proposition 3.2]{_Shah:uniformly_},
cf. \cite[Proposition 3.3.7]{_Kleinbock_etc:Handbook_})

Let $G$ be
the identity component of
a linear algebraic $\Q$-group
$G'$, and
$\Gamma\subset G$ an arithmetic Lie lattice. Let $H\subset G$ be a closed Lie
subgroup generated by
algebraic unipotent
one-parameter subgroups of $G$ contained in $H$.
Let $x:=eH\in G/H$,
where $e\in G$ is the identity of $G$.

Then the closure of the orbit $\Gamma \cdot x$ in $G/H$
is $\Gamma S\cdot x$, where $S\subset G$ is
the identity component of
the smallest
algebraic $\Q$-subgroup of
$G'$
containing $H$.
\endproof

\hfill


\hfill


\remark \label{_Borel_ChH_Remark_}

This result is stated in \cite[Proposition 3.3.7]{_Kleinbock_etc:Handbook_} under the
assumption that $G$ has no $\Q$-characters. In fact, in
view of Borel and Harish-Chandra theorem
(\ref{_Borel-Harish-Chandra_Claim_}), this assumption is redundant.


\hfill


\section{The $\T^4$ and K3 case -- application of Shah's version of Ratner's orbit closure theorem}
\label{_Ratner-K3-case_Section_}

The goal of this section
is to prove
\ref{_dense-orbit-of-sympl-form-K3-case_Theorem_}.

In this section we will assume that $\Lambda$ is
an even unimodular lattice of signature $(p,q)$, $p,q\geq 3$.
Let $V:=\Lambda\otimes_\Z \R$.
The symmetric bilinear form on $\Lambda$ and $V$ will be denoted as before by $(\cdot,\cdot)$.

The group $SO (V)\subset SL(V)$ is a linear algebraic $\Q$-group and $SO^+(V)$ is its identity component.
(Note, however, that $SO^+(V)$ itself is not an algebraic subgroup --
its Zariski closure is $SO(V)$ \cite[p.3]{_Satake_symme_}).

Let $u\in\Lambda$ be a primitive isotropic non-zero vector.

Denote the stabilizers of $u$ in $SO (V)$ and $SO^+(V)$ respectively by $\Gu'$ and $\Gu$:
\[
\Gu' := \{\ g\in SO (V)\ |\ gu=u  \ \}.
\]
\[
\Gu := \{\ g\in SO^+(V)\ |\ gu=u  \ \}.
\]
Each of these sets is the intersection of two Lie subgroups of the group
$SL (V)$: the stabilizer of $u$ in $SL (V)$ and $SO(V)$, or, respectively,
$SO^+(V)$. Since the intersection of two Lie subgroups is
a Lie subgroup, $\Gu'$ is a Lie group and $\Gu$ is its closed Lie subgroup.

Since $u\in\Lambda$, the group $\Gu'$ is a linear algebraic $\Q$-group.

The group $\Gu'$ acts on $u^\bot$ and preserves $u^{\bot,1}$.

In this case
\[
\GLambda = \Gu\cap SO^+ (\Lambda).
\]

Let $y\in u^{\bot,1}$ be a $u$-orthoirrational unit vector.

In order to prove \ref{_dense-orbit-of-sympl-form-K3-case_Theorem_}
we need the following key proposition which will be proved using \ref{_Ratner_arithmetic_Claim_}.


\hfill


\proposition\label{_G-u-Lambda-orbit-dense-in-hyperplane_Proposition_}

The $\GLambda$-orbit of $y$ is dense in $u^{\bot,1}$,
for any $u$-orthoirra\-tio\-nal vector $y\in u^\bot$.


\hfill


\subsection{Density of $\GLambda$-orbits in $u^{\bot,1}$}
\label{_density-of-orbit-K3-case-proof_Subsection_}

\noindent
{\sl Preparations for the proof of \ref{_G-u-Lambda-orbit-dense-in-hyperplane_Proposition_}:}

\bigskip
For the proof of \ref{_G-u-Lambda-orbit-dense-in-hyperplane_Proposition_} we first need a number of preparations.

Set
\[
k:= p+q-1.
\]
We have $\dim_\R u^\bot = k$.

The restriction map $\phi\mapsto \phi|_{u^\bot}$
provides a canonical isomorphism between $\Gu$ and the identity component of the group of isometries of $u^\bot$ fixing $u$.
Indeed, any $(k-1)$-dimensional vector subspace $Z$ of $u^\bot$ transversal to $\Span_\R \{ u\}$ determines a unique vector $v_Z$ such that $v_Z\bot Z$, $(v_Z,v_Z)=1$, $(v_Z,u)=1$. It is easy to see that $v_Z\notin u^\bot$. An isometry of $u^\bot$ fixing $u$ sends $Z$ to another $(k-1)$-dimensional vector subspace $Z'$ of $u^\bot$ transversal to $\Span_\R \{ u\}$ and extends uniquely to an element of $\Gu$ sending $v_Z$ to $v_{Z'}$. This provides an inverse to the restriction map above and shows that it is an isomorphism.
Thus, we can describe elements of $\Gu$ in terms of their restriction to $u^\bot$.

Denote by $\Hy$ the stabilizer of $y$ in $\Gu$:
\[
\Hy := \{\ g\in \Gu\ |\ gy=y\ \}.
\]
This is the intersection of two Lie subgroups of the group $SL (V)$ of isomorphisms of $V$: the stabilizer of $y$ in $SL (V)$ and $\Gu$. Since the intersection of two Lie subgroups is a Lie subgroup, $\Hy$ is a Lie group.

Define
\[
W:= u^\bot/\Span_\R \{u\}.
\]
The bilinear form $(\cdot,\cdot)$ on $V$ induces a non-degenerate bilinear form on $W$ of signature $(p-1,q-1)$. Denote by $SO (W)$ the group of isomorphisms of $W$ preserving the latter bilinear form and by $SO^+ (W)$ its connected component of the identity.

Let $y^\bot_{u^\bot}$ be the orthogonal complement of $y$ in $u^\bot$. Define
\[
W' := y^\bot_{u^\bot} / \Span_\R \{ u\}.
\]
The bilinear form $(\cdot,\cdot)$ on $u^\bot$ induces a non-degenerate bilinear form on $W'$ of signature $(p-2,q-1)$. Denote by $SO (W')$ the group of isomorphisms of $W'$ preserving the latter bilinear form and by $SO^+ (W')$ its connected component of the identity.

We will consider bases $\cB=\{ w_1,\ldots,w_{k-1},u\}$ of $u^\bot$ with $u$ being the last vector of the basis. Such a basis will be called {\bf adapted}.

Since the restriction of $(\cdot,\cdot)$ to $\Span_\R \{ w_1,\ldots,w_{k-1} \}$ is a non-degenerate bilinear form of signature $(p-1,q-1)$, the vectors $w_1,\ldots,w_{k-1}$ in an adapted basis can be assumed to form an orthonormal basis of their span. (A basis is said to be {\bf orthonormal} with respect to an indefinite symmetric bilinear form if it is orthogonal and the square of each basic vector is $\pm 1$).
Such an adapted basis will be called an {\bf adapted orthonormal basis}.

On the other hand, part (B1) of \ref{_unimodular-lattices-divisor-of-any-integral-vector-is-one_Proposition_}
implies that $\Lambda\cap u^\bot$ is a lattice in $u^\bot$ and therefore one can choose an adapted basis $\cB'=\{ w_1,\ldots,w_{k-1},u\}$ so that the vectors $w_1,\ldots,w_{k-1}$ all lie in $\Lambda\cap u^\bot$. Such a basis will be called an {\bf adapted integral basis}.

\bigskip
{\bf Further on let $\cB$ denote an adapted orthonormal basis $\cB$ of $u^\bot$ of the form $\cB=\{ w_1,\ldots,w_{k-2},y,u\}$, with $w_{k-1} =y$.}
\bigskip

There is a natural homomorphism $\Phi: \Gu\to SO^+ (W)$. It is surjective and admits a right inverse. Indeed, any element $\psi\in SO^+ (W)$ induces a unique isometry of the $\Span_\R \{ w_1,\ldots,w_{k-2},y\}$ that extends uniquely to an isometry of $u^\bot$ fixing $u$. Denote this isometry by $\Psi_\cB (\psi)$. One easily checks that such a $\Psi_\cB: SO^+ (W)\to \Gu$ is a right inverse of $\Phi$.

The kernel of $\Phi$ is a Lie subgroup of $\Gu$ that will be denoted by $U$:
\[
U := \ker \Phi.
\]

Clearly,
\begin{equation}
\label{eqn-U-SOplusW-generate-G}
\left.
\begin{aligned}
& \textrm{Any}\ g\in \Gu\ \textrm{can be written in a unique way as}\\
& g=\phi \psi,\  \phi\in \Psi_\cB \left(SO^+ (W)\right), \psi\in U.
\end{aligned}
\right.
\end{equation}

There is also a canonical homomorphism
$\Phi': \Hy\to SO^+ (W')$. It is surjective as well, and admits a right inverse.
Indeed, any element $\psi'\in SO^+ (W')$ induces a unique isometry of the $\Span_\R \{ w_1,\ldots,w_{k-2}\}$ that extends uniquely to an isometry of $u^\bot$ fixing $u$ and $y$. Denote this isometry by $\Psi'_\cB (\psi')$. One easily checks that $\Psi'_\cB: SO^+ (W')\to \Hy$ is a right inverse of $\Phi'$.

The kernel of the homomorphism $\Phi': \Hy\to SO^+ (W')$ is a Lie subgroup of $\Hy$ that will be denoted by $U'$:
\[
U' := \ker \Phi'.
\]
Clearly,
\begin{equation}
\label{eqn-U-prime-SOplusW-prime-generate-H}
\left.
\begin{aligned}
& \textrm{Any}\ h\in \Hy\ \textrm{can be written in a unique way as}\\
& h=\phi' \psi',\  \phi'\in \Psi'_\cB \left(SO^+ (W')\right), \psi'\in U'.
\end{aligned}
\right.
\end{equation}

The elements of $\Gu$ are represented, with respect to $\cB$, by $k\times k$-matrices of the form
\begin{equation}\label{_St_L_block_Equation_}
\begin{pmatrix}
A & \rvline & \begin{matrix} 0 \\ \vdots\\ 0 \end{matrix}\\
\hline
* \dots * &\rvline &  1
\end{pmatrix}
\end{equation}
(Here and further on the asterisks are arbitrary real numbers).
The upper-left $(k-1)\times (k-1)$-block $A$ is a matrix in $SO^+ (p-1,q-1)$.

The elements of $U$ are represented, with respect to $\cB$, by matrices of the form:
\begin{equation}\label{_St_L_block_Equation1_}
\begin{pmatrix}
I_{k-1} & \rvline & \begin{matrix} 0 \\ \vdots\\ 0 \end{matrix}\\
\hline
* \dots * &\rvline &  1
\end{pmatrix}
\end{equation}

The elements of the Lie algebra of $U$ can be then identified with the space of matrices of the form
\[
\begin{pmatrix}
\bigzero & \rvline & \begin{matrix} 0 \\ \vdots\\ 0 \end{matrix}\\
\hline
* \dots * &\rvline &  0
\end{pmatrix}
\]

The elements of $\Hy$ are represented, with respect to $\cB$, by matrices of the form
\begin{equation}\label{_St_L_block_Equation3_}
\begin{pmatrix}
A' & \rvline & \begin{matrix} 0 \\ \vdots\\ 0 \end{matrix} & \rvline & \begin{matrix} 0 \\ \vdots\\ 0 \end{matrix}\\
\hline
0 \dots 0 &\rvline &  1 &\rvline & 0\\
\hline
* \dots * &\rvline &  0 &\rvline & 1
\end{pmatrix}
\end{equation}
The upper-left $(k-2)\times (k-2)$-block $A'$ is a matrix in $SO^+ (p-2,q-1)$.

The elements of $U'$ are represented, with respect to $\cB$, by matrices of the form:
\[
\begin{pmatrix}
I_{k-2} & \rvline & \begin{matrix} 0 \\ \vdots\\ 0 \end{matrix} & \rvline & \begin{matrix} 0 \\ \vdots\\ 0 \end{matrix}\\
\hline
0 \dots 0 &\rvline &  1 &\rvline & 0\\
\hline
* \dots * &\rvline &  0 &\rvline & 1
\end{pmatrix}
\]

The elements of the Lie algebra of $U'$ can be then identified with the space of matrices of the form
\begin{equation}\label{_St_L_block_Equation5_}
\begin{pmatrix}
\bigzero & \rvline & \begin{matrix} 0 \\ \vdots\\ 0 \end{matrix} & \rvline & \begin{matrix} 0 \\ \vdots\\ 0 \end{matrix}\\
\hline
0 \dots 0 &\rvline &  0 &\rvline & 0\\
\hline
* \dots * &\rvline &  0 &\rvline & 0
\end{pmatrix}
\end{equation}


\hfill


\remark

The group $G$ is a (maximal) parabolic subgroup of $SO^+ (V)$ and the decomposition
$G= U\rtimes SO^+ (W)$ is the Levi decomposition of $G$.


\hfill


\bigskip
\noindent
{\sl Lemmas needed for the proof of \ref{_G-u-Lambda-orbit-dense-in-hyperplane_Proposition_}:}

\bigskip

Now we state and prove a number of lemmas that will allow us to verify that
\ref{_Ratner_arithmetic_Claim_} is applicable in our setting and then to apply it.


\hfill


\lemma\label{_Gu-identity-component-of-Gu-prime_Lemma_}

The group $\Gu$ is the identity component of $\Gu'$.


\hfill


\noindent
{\bf Proof of \ref{_Gu-identity-component-of-Gu-prime_Lemma_}.}

It follows from
\eqref{_St_L_block_Equation_}
that
$\Gu$ is connected. Since $\Gu$ is the intersection of $\Gu'\subset SO (V)$
with the identity component $SO^+ (V)$ of $SO (V)$, this implies that $\Gu$ is the identity component of $\Gu'$.
\endproof


\hfill


\lemma\label{_Hy-genertd-by-unipotents_Lemma_}

The group $\Hy$ is generated by
its algebraic
 one-parameter
 unipotent subgroups.


\hfill


\noindent
{\bf Proof of \ref{_Hy-genertd-by-unipotents_Lemma_}.}

The matrices of the form \eqref{_St_L_block_Equation5_} are nilpotent and therefore the corresponding elements of the Lie algebra of $U'$ are $ad$-nilpotent. Moreover, it is easy to see that each element of $U'$ lies in an algebraic
 one-parameter
 unipotent subgroup of $U'$.

By \ref{_SO-p-q-generated-by-unipotents_Proposition_}, since $p,q\geq 3$, the group $SO^+ (p-2,q-1)$ is generated by its algebraic one-parameter
 unipotent subgroups.
Since the isomorphism $SO^+ (W')\cong SO^+ (p-2,q-1)$ is a polynomial map, it sends algebraic
 one-parameter
 unipotent subgroups into algebraic
 one-parameter
 unipotent subgroups. Hence, $SO^+ (W')$ is generated by its algebraic
 one-parameter
 unipotent subgroups. Since $U'$ and $SO^+ (W')$ are generated by their algebraic
 one-parameter
 unipotent subgroups, so is $\Hy$, because of \eqref{eqn-U-prime-SOplusW-prime-generate-H}.

This proves the lemma.
\endproof


\hfill


\lemma\label{_Ky-existence_Lemma_}

There is a unique connected Lie subgroup $\Ky\subset \Gu$ satisfying
\[
\Hy\subsetneqq \Ky\subsetneqq \Gu.
\]
It is formed by all the elements of $\Gu$ preserving both $y^\bot_{u^\bot}$ and send $y$ to vectors of the form $y+cu$, $c\in\R$.


\hfill


\noindent
{\bf Proof of \ref{_Ky-existence_Lemma_}.}

Let us first prove that if a connected Lie group $\Ky$ satisfies $\Hy\subsetneqq \Ky\subsetneqq \Gu$, then it is formed by all the elements of $\Gu$ that preserve both $y^\bot_{u^\bot}$ and $\Span_\R \{ y,u\}$.

Indeed, assume that $\Ky$ is a connected Lie group such that $\Hy\subsetneqq \Ky\subsetneqq \Gu$.
As before, we work with the adapted orthonormal basis $\cB$ and use the identifications
$W=\Span_\R \{ w_1,\ldots,w_{k-2},y \}$
and $W'=\Span_\R \{ w_1,\ldots,w_{k-2} \}$ -- these identifications preserve the bilinear forms and induce isomorphisms
\[
\varsigma_\cB: SO^+ (W)\to SO^+ (p-1,q-1),\ \varsigma'_\cB: SO^+ (W')\to SO^+ (p-2,q-1).
\]
In terms of the matrix representations \eqref{_St_L_block_Equation_}, \eqref{_St_L_block_Equation3_} and the identifications $\varsigma_\cB, \varsigma'_\cB$
above, $\Phi$ can be viewed as a surjective homomorphism $\Phi: \Gu\to SO^+ (p-1,q-1)$
sending each matrix \eqref{_St_L_block_Equation_} to its upper-left $(k-1)\times (k-1)$-block, while $\Phi'$ can be viewed as a surjective homomorphism $\Phi': \Hy\to SO^+ (p-2,q-1)\subset SO^+ (p-1,q-1)$
sending each matrix \eqref{_St_L_block_Equation3_} to its upper-left $(k-2)\times (k-2)$-block, so that $\Phi' = \Phi|_{\Hy}$.

Consequently, $SO^+ (p-2,q-1)\subset \Phi (\Ky)\subset SO^+ (p-1,q-1)$.
Since $\Ky$ is connected, so is $\Phi (\Ky)$. By \cite[Lem. 9.9]{_EV-JTA_}, this implies that either $\Phi (\Ky) = SO^+ (p-2,q-1)$ or $\Phi (\Ky) = SO^+ (p-1,q-1)$.

If $\Phi (\Ky) = SO^+ (p-1,q-1)$, then, in view of \eqref{eqn-U-SOplusW-generate-G} and $\Hy\subsetneqq \Ky\subsetneqq \Gu$,
the matrices representing the elements of $\Ky$ with respect to $\cB$ have to be exactly the matrices of the form
\[
\begin{pmatrix}
A & \rvline & \begin{matrix} 0 \\ \vdots\\ 0 \end{matrix}\\
\hline
* \dots * 0 &\rvline &  1
\end{pmatrix}
\]
where $A\in SO^+ (p-1,q-1)$. However, this collection
 of
 matrices is not a group. This means that $\Phi (\Ky) \neq SO^+ (p-1,q-1)$.

Thus, $\Phi (\Ky) = SO^+ (p-2,q-1)$. In this case, in view of
\eqref{eqn-U-SOplusW-generate-G} and $\Hy\subsetneqq \Ky\subsetneqq \Gu$,
we get that the matrices representing the elements of $\Ky$ with respect to $\cB$ are exactly the matrices of the form
\[
\begin{pmatrix}
A' & \rvline & \begin{matrix} 0 \\ \vdots\\ 0 \end{matrix} & \rvline & \begin{matrix} 0 \\ \vdots\\ 0 \end{matrix}\\
\hline
0 \dots 0 &\rvline &  1 &\rvline & 0\\
\hline
* \dots * &\rvline &  * &\rvline & 1
\end{pmatrix}
\]
where $A'\in SO^+ (p-2,q-1)$.
Such a collection of matrices does form a connected Lie group. It is exactly the group formed by all the elements of $\Gu$ that preserve both $y^\bot_{u^\bot} = \Span\{ w_1,\ldots, w_{k-2},u\}$ and
 sending
 $y$ to vectors of the form $y+cu$, $c\in\R$.

This finishes the proof of the lemma.
\endproof


\hfill


\lemma\label{_Gu-transitivity_Lemma_}

The group $\Gu$ acts transitively on the set of vectors of $u^\bot$ of the same positive length (recall that the length is measured with respect to the bilinear form $(\cdot,\cdot)$).


\hfill


\noindent
{\bf Proof of \ref{_Gu-transitivity_Lemma_}.}

Assume $y'\in u^\bot$, $(y',y')=(y,y)=:r>0$. We want to show that there is an element of $\Gu$ mapping $y$ to $y'$.

Consider a subspace $Z\subset u^\bot$, $\dim_\R Z = \dim_\R u^\bot - 1=k-1$, that is transversal to $u$ and contains $y,y'$. The restriction of $(\cdot,\cdot)$ to $Z$ is a non-degenerate bilinear form of signature $(p-1,q-1)$. The group $SO(Z)$ of isometries of $Z$ is isomorphic to $SO (p-1,q-1)$. In particular, this implies that the identity component $SO^+ (Z)$ of $SO(Z)$ acts transitively on the set of vectors of length $r$ in $Z$. An isometry of $Z$ lying in $SO^+ (Z)$ and sending $y$ to $y'$ can be extended to an isometry of $u^\bot$ fixing $u$ -- that is, to an element of $\Gu$.

This proves the lemma.
\endproof


\hfill


{\bf Assume now that $\cB'=\{ w_1,\ldots,w_{k-1},u\}$ is an adapted integral basis of $u^\bot$.}


\hfill


\lemma\label{_Gu-no-Q-characters_Lemma_}

The group $\Gu$ has no non-trivial $\Q$-characters.


\hfill


\noindent
{\bf Proof of \ref{_Gu-no-Q-characters_Lemma_}.}

Let $\rho: U\to \R^{k-1}$ be an isomorphism
sending each element of $U$, represented, with respect to $\cB'$, by a matrix
\[
\begin{pmatrix}
I_{k-1} & \rvline & \begin{matrix} 0 \\ \vdots\\ 0 \end{matrix}\\
\hline
a_1 \dots a_{k-1} &\rvline &  1
\end{pmatrix}
\]
to $(a_1,\ldots,a_{k-1})$.

Assume that $\chi$ is a $\Q$-character on $\Gu$. Then $\chi|_{\Psi \left(SO^+ (W)\right)}\equiv 1$ -- indeed, $SO^+ (p,q)$ is linear, connected and semi-simple, hence coincides with its commutator subgroup and therefore does not admit any non-trivial characters.
At the same time $\chi|_U$ is a $\Q$-character on $U$ and consequently $f:=\chi|_U \circ \rho^{-1}$ is a $\Q$-character on $\R^{k-1}$ -- that is, $f$ is a polynomial with rational coefficients on $\R^{k-1}$, with values in $\R_{>0}$, satisfying $f(a+b) = f(a)f(b)$ for all $a,b\in\R^{k-1}$. This readily implies that  $f\equiv 1$ and consequently $\chi|_U\equiv 1$. Since $\chi|_U\equiv 1$ and $\chi|_{\psi \left(SO^+ (W)\right)}\equiv 1$, we get
$\chi\equiv 1$ because of \eqref{eqn-U-SOplusW-generate-G}. Thus, $\Gu$ does not admit non-trivial $\Q$-characters.

This proves the lemma.
\endproof


\hfill


\lemma\label{_GLambda-arithmetic-lattice-K3-case_Lemma_}

The group $\GLambda$ is an arithmetic Lie lattice in $\Gu$.


\hfill


\noindent
{\bf Proof of \ref{_GLambda-arithmetic-lattice-K3-case_Lemma_}.}

Follows immediately from \ref{_Gu-identity-component-of-Gu-prime_Lemma_}, \ref{_Gu-no-Q-characters_Lemma_}
and
\ref{_lattice-identity-component_Proposition_}.
\endproof


\hfill


\lemma\label{_H-y-K-y-not-Q-subgps_Lemma_}

Consider the group $\Hy$, written in matrix form in
\eqref{_St_L_block_Equation3_}, and $\Ky$,
defined in \ref{_Ky-existence_Lemma_}. Then $\Hy$ and $\Ky$
are not identity components of
$\Q$-subgroups of
$\Gu'$.


\hfill


\noindent
{\bf Proof of \ref{_H-y-K-y-not-Q-subgps_Lemma_}.}

Set
\[
n:=p+q = k+1.
\]
Without loss of generality, we may assume that $V=\R^n$
and $\Lambda=\Z^n$, $u=(1,0,\ldots,0)\in\Z^n$, and
$(\cdot,\cdot)$ is a bilinear symmetric form of signature
$(p,q)$ on $\R^n$ with an invertible integral matrix.

Clearly, $\Gu'$ is a $\Q$-subgroup.

Let us prove the claim about $\Hy$. Assume, by
contradiction, that $\Hy$ is the identity component of a
$\Q$-subgroup $S$ of
$\Gu'$.

Then the Lie algebras $\Lie (S)$ and $\Lie (\Hy)$ of $S$
and $\Hy$ coincide: $\Lie (S) = \Lie (\Hy)$. Since $S$ is
the common zero level set of a finite number of
polynomials (in the entries of matrices in
$\Gu'$) with
rational coefficients, its Lie algebra $\Lie (S)$ (that
is, the tangent space of $S$ at the identity) is the
common zero level set of a finite number of linear
functions with rational coefficients, hence a Lie
subalgebra of the Lie algebra
$\Lie (\Gu')$ of $\Gu'$
admitting
a basis $A_1,\ldots,A_k$ formed by matrices with rational
coefficients.

Since $\Hy$ is the stabilizer of $u$ and $y$,
the Lie algebra $\goth h$ of $\Hy$ is the space of all
matrices $h\in \goth{so}(V)$ such that $h(y)=h(u)=0$.
This gives $\goth h(V) = \langle u, y\rangle^\bot$.
If $\goth h$ is rational, its image $\goth h(V)$
is also rational, implying that $\langle u, y\rangle^\bot$
is also rational. This is impossible, because
$y$ is $u$-orthoirrational
(\ref{_isotropically-irrational-vector_Definition_}).
This proves that $\Hy$ is not the identity
component of a $\Q$-subgroup.

Let us prove the claim in the case of $\Ky$. Assume, by
contradiction, that $\Ky$ is the identity component of a
$\Q$-subgroup $S$ of
$\Gu'$.
Similarly to the previous case we get that the Lie algebra
$\Lie (S) = \Lie (\Ky)$ admits a basis $A_1,\ldots,A_k$
formed by matrices with rational coefficients.
By definition, $K$ is the subgroup of
$\Gu'$
preserving $u$ and mapping $y$ to $y+cu$, for $c\in \R$.
The Lie algebra $\goth k$ of $K$ is an algebra of all
$h\in \goth{so}(V)$ such that $h(u)=0$ and $h(y)=c u$.
Therefore,
$\goth k$
surjectively maps $V$ to
$y^\bot$; it cannot be rational, because $y^\bot$ is irrational.

This finishes the proof in the case of $\Ky$, as well as
the proof of the lemma.
\endproof


\hfill


\bigskip
\noindent
{\sl Final steps:}

\bigskip

Finally, we are ready to prove
\ref{_G-u-Lambda-orbit-dense-in-hyperplane_Proposition_}.


\hfill


\noindent
{\bf Proof of \ref{_G-u-Lambda-orbit-dense-in-hyperplane_Proposition_}.}

We would like to deduce the wanted result from \ref{_Ratner_arithmetic_Claim_}. Let us check that \ref{_Ratner_arithmetic_Claim_} is applicable here.

Indeed, $\Gu'$ is a linear algebraic $\Q$-group and $G$ is its identity component -- see \ref{_Gu-identity-component-of-Gu-prime_Lemma_}.

The group $\GLambda$ is an arithmetic Lie lattice in $\Gu$ -- see \ref{_GLambda-arithmetic-lattice-K3-case_Lemma_}.

The Lie subgroup $\Hy\subset \Gu$ is closed. It is generated by its algebraic
 one-parameter
 unipotent subgroups -- see \ref{_Hy-genertd-by-unipotents_Lemma_}.

Thus, \ref{_Ratner_arithmetic_Claim_} is applicable and yields that
the closure of the orbit $\GLambda \cdot y$ in $\Gu/\Hy$
is $\GLambda S\cdot y$, where $S\subset \Gu$ is the identity component of the smallest
algebraic $\Q$-subgroup of $\Gu'$ containing $\Hy$.

By \ref{_Ky-existence_Lemma_}, the identity component of $S$ is either $\Hy$, or $\Ky$, or $\Gu$ (the latter case, of course, means that $S=\Gu$).
Since $y\in u^\bot$ is $u$-orthoirrational, \ref{_H-y-K-y-not-Q-subgps_Lemma_} implies that the groups $\Hy$ and $\Ky$ are not identity components of $\Q$-subgroups of $\Gu'$.

Thus, $S=\Gu$ and the closure of the orbit $\GLambda \cdot y$ in $\Gu/\Hy$
is
\[
\GLambda \Gu\cdot y = \Gu y = \Gu/\Hy.
\]
In other words, the orbit $\GLambda \cdot y$ is dense in $\Gu/\Hy$.

On the other hand, by \ref{_Gu-transitivity_Lemma_}, $\Gu$ acts transitively on $u^{\bot,1}$ and therefore $\Gu/\Hy=u^{\bot,1}$.
Thus, the orbit $\GLambda \cdot y$ is dense in $u^{\bot,1}$, which finishes the proof of the proposition.
\endproof


\hfill


\subsection{The density of $SO^+ (\Lambda)$-orbits}
\label{_density-of-SO-plus-Lambda-orbits_Subsection_}

We are now ready to deduce the following corollary of \ref{_G-u-Lambda-orbit-dense-in-hyperplane_Proposition_}.


\hfill


\theorem
\label{_SO-Lambda-orbit-dense-in-hyperplane_Theorem_}

Assume that $\Lambda$ is an even unimodular lattice of signature $(p,q)$, $p,q\geq 3$, which is isotropically transitive.
Let $v\in V$ be an orthoisotropically irrational unit vector.

Then for any non-zero isotropic vector $u\in\Lambda$ the intersection of the $SO^+ (\Lambda)$-orbit of $v$ with $u^{\bot,1}$ is dense in $u^{\bot,1}$.


\hfill


\noindent
{\bf Proof of \ref{_SO-Lambda-orbit-dense-in-hyperplane_Theorem_}.}

Since $v\in V$ is orthoisotropically irrational, there exists a primitive isotropic $x\in\Lambda$ such that $v\in x^\bot$ and $v$ is $x$-orthoirrational.

Let $u\in\Lambda$ be a non-zero isotropic vector. Without loss of generality, we may assume that $u$ is primitive.
Since $\Lambda$ is isotropically transitive, there exists $\phi\in SO^+ (\Lambda)$ mapping $x$ into $u$ and, accordingly, $x^\bot$ into $u^\bot$. The vector $\phi (v)$ then lies in $u^{\bot,1}$ and is $u$-orthoirrational. Applying \ref{_G-u-Lambda-orbit-dense-in-hyperplane_Proposition_} to $\phi (v)\in u^{\bot,1}$ we obtain that the needed result.
\endproof


\hfill


\subsection{Density of the $\Diff^+ (M)$-orbits in cohomology}
\label{_density-of-Diff-orbits-in-cohomology_Subsection_}

Let $M$ be either $\T^4$ or a smooth manifold underlying a complex K3 surface.

The goal of this section
is to formulate
certain density properties of $\Diff^+ (M)$-orbits in $H^2 (M;\R)$. We will then use these properties in Section~\ref{_dense-orbits-of-forms-pfs_Subsection_} to prove 
\ref{_dense-orbit-of-sympl-form-K3-case_Theorem_} about density properties of
the corresponding $\Diff^+ (M)$-orbits in $\SymplKone (M)$.


\hfill


\theorem\label{_dense-orbit-of-cohom-class-K3-case_Proposition_}

Assume
$y\in P = \{\ x\in H^2 (M;\R)\ |\ (x,x)=1\ \}$ is orthoisotropically irrational.

Then for any isotropic $u\in H^2 (M;\Z)$ the intersection of the $\Diff^+ (M)$-orbit of $y$ with $u^\bot\cap P$ is dense in $u^\bot\cap P$. In particular, it is dense in $P$ (for $u=0$).


\hfill


\noindent
{\bf Proof of \ref{_dense-orbit-of-cohom-class-K3-case_Proposition_}.}

We are going to reduce the claim to \ref{_SO-Lambda-orbit-dense-in-hyperplane_Theorem_}.

In the case of $M=\T^4$ we identify $H^2(\T^4; \R)$, with the intersection product on it, with the quadratic space $V:= \big(\bigwedge^2
\R^4\big)^*$ of bilinear forms on $\R^4$ (with the appropriate bilinear symmetric form on $V$), and $H^2(\T^4; \Z)$ with the appropriate quadratic lattice $\Lambda\subset V$. Since the $\Diff^+ (\T^4)$-action on $H^2(\T^4; \R)$ preserves $H^2 (\T^4; \Z)$, it induces a homomorphism from $\Diff^+ (\T^4)$ to the group $SO^+ (\Lambda)$. The image of this homomorphism can be identified with $SL(4,\Z)$ that acts naturally on $V= \big(\bigwedge^2
\R^4\big)^*$, as the restriction of the natural $SL (4,\R)$-action on
$V=\big(\bigwedge^2
\R^4\big)^*$.
The quadratic lattice $\Lambda=H^2(\T^4; \R)$ is
isomorphic to $\bfU\oplus \bfU\oplus\bfU$ -- in
particular, it is even, unimodular and of signature
$(3,3)$. It is also isotropically transitive, by
\ref{_U-U-isotropically transitive_Example_}.

The group $SO^+ (V)$ is then isomorphic to
$SL (4,\R)/\pm 1$ and the isomorphism sends
$SO^+ (\Lambda)\subset SO^+ (V)$ to $SL (4,\Z)$.
The
orbits of the $\Diff^+ (\T^4)$-action on $H^2(\T^4; \R)$ are identified with the orbits of the action of $SO^+ (\Lambda)=SL (4,\Z)$ on $V=\big(\bigwedge^2
\R^4\big)^*$.

In the case where $M$ is a smooth manifold underlying a complex K3 surface
the group $\Lambda= H^2 (M;\Z)$, equipped with the intersection form, is a quadratic lattice isomorphic to $\bfU \oplus \bfU \oplus \bfU \oplus
(-E_8)\oplus (-E_8)$ (see e.g. \cite{_Geom-K3-Asterisque_} for the definition of the integral lattice $E_8$ and for the proof). In particular, it is
an even unimodular quadratic lattice of signature $(3,19)$ in $V := H^2 (M;\R)$. It is also isotropically transitive, by \ref{_U-U-isotropically transitive_Example_}.

Since the $\Diff^+ (M)$-action on $H^2(M;\R)$ preserves $H^2(M;\Z)$ and the intersection form, it induces a homomorphism $\Diff^+ (M)\to SO (\Lambda)$. The image of this homomorphism is
$SO^+ (\Lambda)$ -- see \cite{_Borcea-MathAnn86_}, \cite{_Donaldson-Topology90_}. Thus, the
orbits of the $\Diff^+ (M)$-action on $H^2(M; \R)$ are the orbits of the natural action of $SO^+ (\Lambda)$ on $V$.

Summing up, we see that if $M=\T^4$ or a smooth manifold underlying a complex K3 surface, it is enough to prove that
for any isotropic $u\in H^2 (M;\Z)$ the intersection of the $SO^+ (\Lambda)$-orbit of $y$ with $u^\bot\cap P$ is dense in $u^\bot\cap P$.

In the case $u=0$ the density of the $\Diff^+ (M)$-orbit of $y$ in $u^\bot\cap P=P$ was previously proved in \cite[Thm. 9.8]{_EV-JTA_}. (Note that since the unit vector $y$ is orthoisotropically irrational, it is not proportional to an integral cohomology class, and thus the result of \cite[Thm. 9.8]{_EV-JTA_} applies).

Therefore, we only need to prove the theorem in the case $u\neq 0$. In this case, the result follows immediately from the discussion above and \ref{_SO-Lambda-orbit-dense-in-hyperplane_Theorem_}.
\endproof


\hfill


\section{K\"ahler-type symplectic forms on K3 surfaces}
\label{_K3-surfaces_Section_}

In this section we discuss K\"ahler-type symplectic forms on K3 surfaces.
This will be used further on in Section~\ref{_dense-orbits-of-forms-pfs_Subsection_} to prove \ref{_dense-orbit-of-sympl-form-K3-case_Theorem_}.

Let $M$ be a K3 surface equipped with the standard orientation.

Denote by
$Q: H^2 (M;\C)\times H^2 (M;\C)\to\C$ the Hermitian intersection form:
\[
Q (a,b):= \int_M a\cup \bar{b}.
\]

Recall that a {\bf hyperk\"ahler structure} on $M$ is a collection
\[
(\omega_1,\omega_2,\omega_3, I_1,I_2,I_3)
\]
of three complex structures $I_1,I_2,I_3$, compatible with the orientation and satisfying the
quaternionic relations, and three symplectic forms
$\omega_1,\omega_2,\omega_3$, compatible, respectively, with
$I_1,I_2,I_3$, so that the three Riemannian metrics $\omega_i
(\cdot, I_i\cdot)$, $i=1,2,3$, coincide.
A symplectic form, or a complex structure, on $M$ is said to be of
{\bf hyperk\"ahler type}, if it appears in {\it some} hyperk\"ahler
structure.


\hfill


\proposition
\label{_Kahler-type-on-K3-is-hyperkahler-type_Proposition_}\\
Let $M$ be a K3 surface equipped with the standard orientation.

Then:\\
 - Any complex structure on $M$ (compatible with the orientation)
is of hyperk\"ahler type (and, in particular, of K\"ahler type).\\
- A symplectic form on $M$
 (compatible with the orientation)
 is of K\"ahler-type if and only if
it is of hyperk\"ahler type.


\hfill


\noindent
{\bf Proof of \ref{_Kahler-type-on-K3-is-hyperkahler-type_Proposition_}.}

Let $J$ be a complex structure on $M$.

By a theorem of Siu \cite{_Siu_} (that corrected a
previous result of Todorov \cite{_Todorov_}), $J$ is of K\"ahler type.
(Alternatively, this can be deduced from later results of Buchdahl and Lamari -- see \cite{_Buchdahl:surfaces_}, \cite{_Lamari_} -- who showed
that any complex surface with an even first Betti number admits a K\"ahler structure).

A theorem of Friedman-Morgan \cite[p.495, Cor. 3.5]{_Friedman-Morgan-book_}, based on Donaldson's results on gauge theory, says that any complex surface diffeomorphic to a K3 surface is itself K3  (alternatively, the latter claim can be deduced from the results of Taubes on Seiberg-Witten invariants -- see e.g. \cite[Example 3.13]{_Salamon-uniqueness_}).
Hence, $c_1 (M,J) = 0$ and $(M,J)$ is a complex K3 surface.

Let $\omega$ be a K\"ahler form on $(M,J)$. Then, by Yau's theorem (formerly the Calabi conjecture) \cite{_Yau_}, since $c_1 (M,J) = 0$, the complex structure $J$ belongs to a K\"ahler structure $(\omega',J)$ on $M$ for which the corresponding Hermitian metric $\omega' (\cdot, J\cdot) + \sqrt{-1} \omega' (\cdot,\cdot)$ is Ricci-flat, and moreover $[\omega']=[\omega]$. Consequently, $J$ and $\omega'$ are of hyperk\"ahler type -- see \cite{_Beauville_}, cf. \cite[Thm. 6.40]{_Besse:Einst_Manifo_}. Since $\omega$ and $\omega'$ are cohomologous symplectic forms compatible with the same $J$, by a theorem of Moser \cite{_Moser_}, they can be identified by an isotopy of $M$. Therefore, since $\omega'$ is of hyperk\"ahler type, so is $\omega$.
\endproof


\hfill


Let $J$ be a complex structure on $M$. By \ref{_Kahler-type-on-K3-is-hyperkahler-type_Proposition_}, $J$ if of K\"ahler type and therefore induces a Hodge decomposition\footnote{Note that although the proof of the
existence of the Hodge decomposition involves the whole K\"ahler
structure of which $J$ is a part, one can show -- see e.g.
\cite[Vol.1, Prop. 6.11]{_Voisin-Hodge_} --
that, in fact, the Hodge decomposition depends only on $J$.} of $H^2 (M;\C)$:
\[
H^2 (M;\C) = H^{2,0} (M,J)\oplus H^{1,1} (M,J)\oplus H^{0,2} (M,J).
\]
Moreover, $\dim_\R H^{2,0} (M,J) = 1$
(see e.g. \cite{_Geom-K3-Asterisque_}) and the space $H^{2,0} (M,J)$ completely determines the full Hodge decomposition of $H^2 (M;\C)$:
\[
H^{0,2} (M,J) = \overline{H^{2,0} (M,J)}
\]
and $H^{1,1} (M,J)$ is the orthogonal complement of $H^{2,0} (M,J)\oplus H^{0,2} (M,J)$ with respect to
$Q$.


\hfill


\section{The $\T^4$ and K3 case -- density of the $\Diff^+ (M)$-orbits in the spa\-ce of forms}
\label{_density-of-orbits-of-forms_Section_}

The goal of this section is to prove
\ref{_dense-orbit-of-sympl-form-K3-case_Theorem_}.


\hfill


\subsection{Symplectic K\"ahler-type Teichmuller space and the action of the mapping class group on it}
\label{_Sympl-Teichm-space-mapp-class-gr_Subsection_}

Let $M^{2n}$ be either $\T^{2n}$ or a smooth manifold underlying a complex K3 surface. Let $\Diff_0 (M)$ be the connected component of the identity in
 the group $\Diff^+ (M)$ of orientation-preserving diffeomorphisms of $M$.

Define the {\bf symplectic K\"ahler-type Teichmuller space} $\TeichsK (M)$ by
\[
\TeichsK (M) := \SymplKone (M)/\Diff_0 (M).
\]
Equip $\TeichsK (M)$ with the quotient topology. For $\omega\in \SymplKone (M)$ denote the corresponding element of $\TeichsK (M)$ by $\{ \omega\}$.

Define {\bf the symplectic period map} $\Pers: \TeichsK (M)\to H^2 (M;\R)$ by
\[
\Pers (\{ \omega\}) := [\omega].
\]
Using Moser's stability theorem for
symplectic structures, it is not hard to obtain
that $\TeichsK (M)$ is a finite-dimensional manifold
(possibly non-Hausdorff)
and $\Pers$
is a local diffeomorphism
(\cite{_Moser_},
\cite[Proposition 3.1]{_Fricke_Habermann:symp-moduli_}).

The image of $\Pers$ is exactly $P=\{ x\in H^2 (M;\R)\ |\ x^n =1 \}$.
Indeed, the inclusion $\Im \Pers\subset P$ is obvious. The
inclusion $P\subset Im \Pers$ in the torus case is easy:
any point in $P$ is the cohomology class of a linear
symplectic form on the torus; such a symplectic form is of
K\"ahler type. In the K3 case the equality $P= Im \Pers$
follows, for instance, from \cite[Thm. 5.1]{_Am-Ver_},
where it is proven for all compact hyperk\"ahler manifolds of
maximal holonomy.

Define {\it the mapping class group of $M$} by
\[
\bGamma := \Diff^+ (M)/\Diff_0 (M).
\]

The group $\Diff^+ (M)$
(and, consequently, the mapping class group $\bGamma$)
acts naturally on $\TeichsK (M)$
and on $H^2 (M;\R)$.
The period map $\Pers$ is equivariant with respect to
these actions.

We will now prove \ref{_DiffH-acts-transitively-on-cohomologous-Kahler-type-forms_Proposition_}. For convenience, we restate it her.


\hfill


\proposition (=\ref{_DiffH-acts-transitively-on-cohomologous-Kahler-type-forms_Proposition_})
\label{_DiffH-acts-transitively-on-cohomologous-Kahler-type-forms-copy_Proposition_}

Let $M^{2n}$ be either $\T^{2n}$ or a smooth manifold underlying a complex K3 surface.

Then any two K\"ahler-type symplectic forms on $M$ (compatible with the orientation of $M$) can be mapped
into each other by a diffeomorphism of $M$ acting trivially on homology.


\hfill


\noindent
{\bf Proof of \ref{_DiffH-acts-transitively-on-cohomologous-Kahler-type-forms-copy_Proposition_} (=\ref{_DiffH-acts-transitively-on-cohomologous-Kahler-type-forms_Proposition_}).}

In the case of $M=T^{2n}$ this is proved in \cite[Erratum]{_EV-JTA_} (it follows from the fact that any K\"ahler-type symplectic form on $\T^{2n}$ can be identified with a linear form by a diffeomorphism {\sl acting as the identity on $H^* (\T^{2n})$}).

Assume that $M$ is a smooth manifold underlying a complex K3 surface.
Let $\omega',\omega''\in\SymplK (M)$ so that $[\omega'] = [\omega'']$.

By \ref{_Kahler-type-on-K3-is-hyperkahler-type_Proposition_}, $\omega'$, $\omega''$ can be included in hyperk\"ahler structures
\[
(\omega'=:\omega'_1,\omega'_2,\omega'_3,I'_1,I'_2,I'_3), (\omega''=:\omega''_1,\omega''_2,\omega''_3,I''_1,I''_2,I''_3).
\]
One can assume without loss of generality that $\Span_\R \{ [\omega'_1], [\omega'_2],[\omega'_3]\} = \Span_\R \{ [\omega''_1], [\omega''_2],[\omega''_3]\}$ -- see \cite[Thm. 4.9 and the proof of Thm. 5.1]{_Am-Ver_}.
Denote
\[
W:= \Span_\R \{ [\omega'_1], [\omega'_2],[\omega'_3]\} = \Span_\R \{ [\omega''_1], [\omega''_2],[\omega''_3]\}.
\]
The orthogonal complement of $[\omega']$ in $W$ determines $H^{2,0} (M,I'_1)$ and, similarly, the orthogonal complement of $[\omega'']$ in $W$ determines $H^{2,0} (M,I''_1)$ -- see  \cite[the proof of Thm. 4.9]{_Am-Ver_}.
Since $[\omega'] = [\omega'']$, we get that
\[
H^{2,0} (M,I'_1) = H^{2,0} (M,I''_1).
\]
In view of the discussion at the end of Section~\ref{_K3-surfaces_Section_}, this means that the Hodge decompositions of $H^2 (M;\C)$ induced by $I'_1$ and $I''_1$ coincide.
Therefore, by the global Torelli theorem for K3 surfaces (see \cite{PS-S}, \cite{_Burns-Rapoport_},
\cite{_Siu_}, cf. \cite[p.96]{_Geom-K3-Asterisque_}), there exists a biholomorphism $(M,I'_1)\to (M,I''_2)$ which acts as identity on $H^* (M)$ and, consequently, preserves the orientation. Such a biholomorphism maps $\omega''_1$ to a symplectic form $\eta$ on $M$ compatible with $I'_1$ and cohomologous to $\omega'_1$. It follows then from the theorem of Moser \cite{_Moser_} that there exists an isotopy of $M$ mapping $\eta$ to $\omega'_1$. Therefore, there exists a diffeomorphism of $M$ acting trivially on the homology and mapping $\omega''=\omega''_1$ to $\omega'=\omega'_1$.
\endproof


\hfill


\proposition\label{_Gamma-acts-transitively-on-conn-comps-of-Teichs_Proposition_}

Let $M^{2n}$ be either $\T^{2n}$ or a smooth manifold underlying a complex K3 surface.

Then the mapping class group $\bGamma$ (and, consequently, $\Diff^+ (M)$) acts transitively on the set of connected components of $\TeichsK (M)$.


\hfill


\noindent
{\bf Proof of \ref{_Gamma-acts-transitively-on-conn-comps-of-Teichs_Proposition_}.}

Let $\TeichsK' (M)$ and $\TeichsK'' (M)$ be two connected components of $\TeichsK (M)$. Since $\Pers: \TeichsK' (M)\to P$ and $\Pers: \TeichsK'' (M)\to P$ are diffeomorphisms, there exist $\omega',\omega''\in\SymplK (M)$ such that $\{ \omega'\}\in \TeichsK' (M)$, $\{ \omega''\}\in \TeichsK'' (M)$, $[\omega'] = [\omega'']$.

By \ref{_DiffH-acts-transitively-on-cohomologous-Kahler-type-forms-copy_Proposition_}, there exists an element of $\Diff^+ (M)$ mapping $\omega''$ to $\omega'$ and hence an element of $\bGamma$ mapping $\{ \omega''\}$ to $\{ \omega'\}$.

This finishes the proof of the proposition.
\endproof

\hfill


\subsection{The proof of
\ref{_dense-orbit-of-sympl-form-K3-case_Theorem_}}
\label{_dense-orbits-of-forms-pfs_Subsection_}


\hfill


Now we can finally prove
\ref{_dense-orbit-of-sympl-form-K3-case_Theorem_}.
 We restate it here for convenience.


\hfill


\theorem (= \ref{_dense-orbit-of-sympl-form-K3-case_Theorem_})
\label{_dense-orbit-of-sympl-form-K3-case-copy_Theorem_}

Assume that $\omega_0\in \SymplKone (M)$ so that the cohomology class $[\omega_0]\in H^2 (M;\R)$ is orthoisotropically irrational.

Then for any isotropic $u\in H^2 (M;\Z)$
the intersection of the $\Diff^+ (M)$-orbit of $\omega_0$ with $\Subot (M)$
is dense in $\Subot (M)$. In particular, it is dense in $\SymplKone (M)$ (for $u=0$).


\hfill


\noindent
{\bf Proof of \ref{_dense-orbit-of-sympl-form-K3-case-copy_Theorem_} (= \ref{_dense-orbit-of-sympl-form-K3-case_Theorem_}).}

The map $\Pers: \TeichsK (M)\to P$ is a surjective local diffeomorphism. Consequently, its restriction to the set $\bigg\{ \{ \omega\}\ |\ \omega\in \Subot\ \bigg\}$ is a local diffeomorphism between this set and $a^\bot\subset P$. The union of the domains and the union of the targets of these local diffeomorphisms for different $a$ are $\Diff^+ (M)$-invariant and the local diffeomorphisms are $\Diff^+ (M)$-equivariant.

Now the theorem
follows from
\ref{_Gamma-acts-transitively-on-conn-comps-of-Teichs_Proposition_}
and
\ref{_dense-orbit-of-cohom-class-K3-case_Proposition_}.
\endproof


\hfill


\section{Proof of the main theorem}
\label{_pfs-main-thm_Section_}

In this section we deduce \ref{_Main_Theorem_} from \ref{_GLambda-orbit-torus-case_Proposition_} and \ref{_dense-orbit-of-sympl-form-K3-case_Theorem_}.

The following proposition is well-known and was already used (without proof) in \cite[Lem. 7.21]{_Sheridan-Smith_}. We present here its proof for completeness.


\hfill


\proposition\label{_Deforming-zero-homologous-Lagr-submfds_Proposition_}

Assume $(M^{2n},\omega_0)$ is a symplectic manifold and $L^n\subset (M^{2n},\omega_0)$ is a closed (that is, compact without boundary) Lagrangian submanifold. As before, let $\Sympl (M)$ be the space of symplectic forms on $M$ equipped with the $C^\infty$-topology.

Then there exists a neighborhood $\cU$ of $\omega_0$ in $\Sympl (M)$ so that for any symplectic form $\omega_1\in\cU$ satisfying $[\omega_1|_L] =0$, the forms $\omega_t:= (1-t)\omega_0 + t\omega_1$, $0\leq t\leq 1$, are symplectic and there exists an isotopy $\{L_t\}_{0\leq t\leq 1}$, $L_0=L$, of $L$ such that each $L_t$ is a Lagrangian submanifold of $(M,\omega_t)$.


\hfill


\noindent
{\bf Proof of \ref{_Deforming-zero-homologous-Lagr-submfds_Proposition_}.}

Fix a Riemannian metric on $M$ and a closed tubular neighborhood $U$ of $L$. Let $\| \cdot\|$ denote the $C^0$-norm of differential forms on $U$ measured with respect to the Riemannian metric.

Let $\omega_1$ be a symplectic form on $M$ such that $[\omega_1|_L] =0$. Then there exists a 1-form $\sigma$ on $U$ such that $d\sigma = (\omega_1 - \omega_0)|_U$. Moreover, one can choose $\sigma$ so that
\begin{equation}
\label{eqn-sigma-norm}
\| \sigma\|\leq C \| (\omega_1 - \omega_0)|_U \|,
\end{equation}
where $C$ depends only on $U$ -- this can be done as in \cite[Lem. 9.3]{_EV-JTA_}.

Let $\cU$ be a neighborhood of $\omega_0$ in $\Sympl (M)$ such that for any $\omega_1\in \cU$ the forms $\omega_t:= (1-t)\omega_0 + t\omega_1 = \omega_0 + t (\omega_1-\omega_0)$, $0\leq t\leq 1$, are symplectic.

Let $\omega_1\in \cU$ and let $X_t$, $0\leq t\leq 1$, be
the vector field on $U$ defined by $\omega_t (X_t, \cdot)
= \sigma (\cdot)$. Shrinking $\cU$, if necessary, to a
smaller neighborhood of $\omega_0$ in $\Sympl (M)$, and
using \eqref{eqn-sigma-norm}, we can assume, without loss
of generality, that the vector field $X_t$ is sufficiently
$C^0$-small. Then its flow $\phi_t$ is defined for
$0\leq t\leq 1$ on a neighborhood of $L$ in $U$.

By a calculation underlying Moser's method \cite{_Moser_}, the flow $\phi_t$, $0\leq t\leq 1$, satisfies (wherever it is defined for all $0\leq t\leq 1$) the condition $\phi_t^* \omega_t = \omega_0$ for all $0\leq t\leq 1$. Thus, $L_t := \phi_t (L)$ is well-defined and Lagrangian with respect to $\omega_t$ for all $0\leq t\leq 1$.
\endproof


\hfill


Further on we will need the following easy lemma.


\hfill

\lemma\label{_primi_subla_transi_Lemma_}

Let $R$ be the set of rational subspaces of $\R^m$ of dimension $k>0$.
Then $SL(m,\Z)$ acts transitively on $R$.


\hfill


\proof

Let $N \subset \Z^m$ be a rational subspace of $\R^m$ of dimension $k>0$.
Then one can choose a basis $x_1,\ldots, x_k$ of the finite-rank free group $N\cap \Z^m$ which is a basis of $N\cap \Q^m$ over $\Q$ (because any vector in $N\cap \Z^m$ has an integral multiple that lies in $N\cap \Q^m$) and hence also a basis of $N$ over $\R$.

Since $SL(m,\Z)$ acts transitively on the set of oriented
bases of $\Z^m$, it would suffice to show that the collection
$x_1,\ldots, x_k$ can be extended to a basis of $\Z^m$.
Since $N$ is primitive,
the quotient $\Z^m/(N\cap \Z^m)$ has no torsion, meaning that it is a finite rank free abelian group.
Choose a basis
$\bar y_1, \ldots, \bar y_{m-k}$  of $\Z^m/(N\cap \Z^m)$ and lift these elements to
$y_1, \ldots, y_{m-k}\in \Z^m$. Then
$x_1, \ldots, x_k, y_1, \ldots, y_{m-k}$ is a basis of $\Z^m$.

This finishes the proof of the lemma.
\endproof


\hfill


Now we are ready to prove \ref{_Main_Theorem_}.
 We restate it here for convenience.


\hfill


\theorem (=\ref{_Main_Theorem_})\label{_Main_Theorem-copy_}

Assume $M$, $\dim_\R M = 2n$, is either an even-dimensional torus or a smooth manifold underlying a K3 surface. Let
$\omega$ be a K\"ahler-type symplectic form on $M$.

Then for any Maslov-zero Lagrangian torus $L\subset (M, \omega)$ the homology class $[L] \in H_n (M;\Z)$ is non-zero
and primitive.


\hfill


\noindent
{\bf Proof of \ref{_Main_Theorem-copy_} (= \ref{_Main_Theorem_})
 -- the torus case}.

Let $x_1,\ldots,x_{2n}$ be the standard coordinates on $\R^{2n}$.

Let
$\omega$ be a K\"ahler-type symplectic form on $\T^{2n}$ and let $L\subset (\T^{2n}, \omega)$
be a Maslov-zero Lagrangian torus $L\subset (\T^{2n}, \omega)$.

We need to prove that the homology class $[L] \in H_n (\T^{2n};\Z)$ is non-zero and primitive.

Since any K\"ahler-type symplectic form on $\T^{2n}$ can be mapped into a linear symplectic form by a diffeomorphism of $\T^{2n}$
(see \cite{_EV-JTA_}) and since the required property of $L$ is invariant under the rescaling of $\omega$ by a constant factor, we may assume without loss of generality that $\omega$ is a linear symplectic form of total volume $1$. Let $\tomega\in\tSlin$ be the lift of $\omega$ to $\R^{2n}$.

Let us identify $\pi_1 (\T^{2n})\cong H_1 (\T^{2n})\cong \Z^{2n}\subset \R^{2n} \cong H_1 (\T^{2n};\R)$ (recall that $\T^{2n} = \R^{2n}/\Z^{2n}$).
Set
\[
Z:=Im\, \left(\pi_1 (L)\to \pi_1 (\T^{2n})\right)\subset\pi_1 (\T^{2n})\cong
\]
\[
\cong
Im\, \left(H_1 (L)\to H_1 (\T^{2n})\right)\subset H_1 (\T^{2n})
\cong
\Z^{2n},
\]
\[
l:= Im\, \left(H_1 (L;\R)\to H_1 (\T^{2n};\R)\right)\subset H_1 (\T^{2n};\R)\cong\R^{2n}.
\]
Let $k:= \dim l$.

Then $Z\subset l\cap \Z^{2n}$ and $l$ is the smallest rational subspace of $\R^{2n}$ containing $Z$. Moreover, $l$ is isotropic with respect to $\tomega$, implying that $0\leq k\leq n$. In fact, $[L]\neq 0$ if and only if $\dim l = n$.

Denote
\[
L_l := \frac{l}{l\cap\Z^{2n}}.
\]
It is a $k$-dimensional torus.

Let $x_1,\ldots,x_{2n}$ be the standard coordinates on
$\R^{2n}$. Since $SL (2n,\Z)$ acts transitively on the set
of integral subspaces of $\R^{2n}$ of dimension $k$
(\ref{_primi_subla_transi_Lemma_}) and maps each linear symplectic form on
$\T^{2n}$ to a linear symplectic form of the same total
volume, we may assume, without loss of generality, that
$l$ is either $\{ 0\}$, if $k=0$, or the
$(x_1,\ldots,x_k)$-coordinate plane, if $1\leq k\leq n$.
Let $l'$ be the $(x_{k+1},\ldots,x_{2n})$-coordinate plane.

Let us show that $[L]\neq 0$. The proof follows the argument in \cite{_Abouzaid-Smith_} for the case of $\T^4$ and the standard Darboux symplectic form on $\T^4$.

Namely, assume, by contradiction that $[L]=0$, or, equivalently, $k < n$.
The natural projection
\[
L_l \times l'\to L_l \times \frac{l'}{l'\cap\Z^{2n}} = \T^{2n}
\]
is a covering. Consider the symplectic form $\homega$ on $L_l \times l'$ which is
the lift of $\omega$. Since
a lift of a Maslov-zero Lagrangian to a covering symplectic space is again a Maslov-zero Lagrangian,
the torus $L$ lifts to a Maslov-zero Lagrangian torus $\hL\subset (L_l \times l',\homega)$. It is also easy to see that the symplectic manifold $(L_l \times l',\homega)$
is convex at infinity in the sense of \cite{_Eliash-Gr-convex_}.

On the other hand, since $l$ is isotropic and $\dim l < n < \dim l'$, we have $l^\bot\cap l'\neq \{ 0\}$,
where $l^\bot$ is the $\omega$-orthogonal complement of $l$. This implies that there exists a Hamiltonian on $L_l \times l'$, which is a linear function of
$x_{k+1},\ldots,x_{2n}$, whose constant Hamiltonian vector field is parallel to $l'$. Consequently, any compact subset of $L_l \times l'$ -- and, in particular, $\hL$ -- is displaceable by a Hamiltonian isotopy. However, Maslov-zero Lagrangian torus cannot be displaceable by a theorem of Fukaya \cite[Thm. 12.2]{_Fukaya-Montreal2004_}.

This leads to a contradiction and shows that $[L]\neq 0$.

Let us now show that $[L]$ is primitive.

As we have shown, $[L]\neq 0$. Consequently, $\pi_1 (L)$ is identified with $Z$ and $l$ and $l'$ are the $(x_1,\ldots,x_n)$ and $(x_{n+1},\ldots,x_{2n})$-coordinate planes. In particular, $l$ and $l'$ are complementary $n$-dimensional rational subspaces. Let $q$ be the index of $Z$ in $l\cap\Z^{2n}$. The homology class $[L]$ is primitive if and only if $q=1$.

Let us first prove the primitiveness of $[L]$ in the case where $\omega$ is $(l,l')$-Lagrangian split. The proof in this case follows the argument in \cite{_Abouzaid-Smith_} for the case of $\T^4$ and the standard Darboux symplectic form on $\T^4$.

Consider the linear Lagrangian torus $L_l := l/(l\cap \Z^{2n})\subset (\T^{2n},\omega)$ and the covering map $L_l\times l'\to L_l\times l'/(l'\cap\Z^{2n}) = \T^{2n}$.
Also consider the $q$-fold covering map $l/Z\to L_l=l/(l\cap \Z^{2n})$ -- its direct product with the identity is a $q$-fold covering map $l/Z\times l'\to L_l\times l'$. Let $\pi: l/Z\times l'\to\T^{2n}$ be the composition of the covering maps $l/Z\times l'\to L_l\times l'$ and
$L_l\times l'\to \T^{2n}$. Let $\homega$ be a symplectic form on $l/Z\times l'$ which is the lift of $\omega$ under $\pi$. It is easy to see that the symplectic manifold $(l/Z\times l',\homega)$
is convex at infinity in the sense of \cite{_Eliash-Gr-convex_}.

There are $q$ disjoint lifts of $L$ to $l/Z\times l'$ and they all are Maslov-zero Lagrangian tori in $(l/Z\times l',\homega)$ mapped into each other by diffeomorphisms of $l/Z\times l'$ induced by parallel translations in $l/Z$. Since $\omega$ is $(l,l')$-Lagrangian split, such parallel translations are Hamiltonian symplectomorphisms of $(l/Z\times l',\homega)$ generated by Hamiltonians that are linear functions of the coordinates $(x_{n+1},\ldots,x_{2n})$ on $l/Z\times l'$. Thus, $q>1$ would have implied that there exists a displaceable Maslov-zero Lagrangian torus
in $(l/Z\times l',\homega)$, in contradiction to the same theorem of Fukaya \cite[Thm. 12.2]{_Fukaya-Montreal2004_} that we have already used above.

This proves the primitiveness of $[L]$ in the case where $\omega$ is $(l,l')$-Lag\-ran\-gian split.

Now let $\omega$ be an arbitrary linear symplectic form on $\T^{2n}$ of total volume $1$.

Assume, by contradiction,
that $[L]$ is not primitive. It follows from \ref{_Deforming-zero-homologous-Lagr-submfds_Proposition_} that for any linear symplectic form
$\omega'$ of total volume $1$, which is sufficiently close to $\omega$ and satisfies $[\omega'|_L]=0$ (for linear symplectic forms $\omega'$ this holds if and only if $\omega'|_l \equiv 0$), the symplectic manifold $(\T^{2n},\omega')$ admits a Maslov-zero Lagrangian torus whose homology class is not primitive. By \ref{_GLambda-orbit-torus-case_Proposition_}, one can find such an $\omega'$ which is the image of an $(l,l')$-Lagrangian split form $\omega''$ under a diffeomorphism of $\T^{2n}$ given by an element of $SL (2n,\Z)$ acting trivially on $l$. Hence, $(\T^{2n}, \omega'')$ too admits a Maslov-zero Lagrangian torus whose homology class is not primitive, in contradiction to the primitiveness in the $(l,l')$-Lagrangian split case that we have proved above.

This finishes the proof of the theorem.
\endproof


\hfill


\noindent
{\bf Proof of \ref{_Main_Theorem-copy_}
 (= \ref{_Main_Theorem_})
 -- the K3 case}.

Let $M$ be a smooth manifold underlying
a smooth K3 surface.

It follows from \cite[Thm. 1.3]{_Sheridan-Smith_} that for
a certain K\"ahler-type complex structure $J$ on $M$ (a so-called ``mirror quartic") there exists a symplectic form $\omega_0\in \SymplKone (M)$, compatible with $J$, such that $(M,\omega)$ does not admit a Maslov-zero Lagrangian torus in the zero or a non-primitive homology class.
In the terminology of \cite{_Sheridan-Smith_}, such an $\omega_0$ can be chosen to be
``ambient irrational". In view of \cite[Example
  3.10]{_Sheridan-Smith_}, the latter property of
$\omega_0$ implies that the quadratic lattice
$[\omega_0]^\bot\cap H^2 (M;\Z)$ is isomorphic to
$\bfU\oplus\bfU\oplus\bfU\oplus \langle 4\rangle$.
(Here
$\langle 4\rangle$ is the quadratic lattice formed by
$\Z$ equipped with the symmetric bilinear form $(\cdot,\cdot)$ defined by $(1,1)=4$).

We claim that $[\omega_0]$ is orthoisotropically irrational.

Indeed, first, $[\omega_0]^\bot$ contains a non-zero primitive isotropic element $u\in H^2 (M;\Z)$ (from a $\bfU$ summand), hence $[\omega_0]\in u^\bot$.

Second, let us show that for each such $u$ we have
\[
\Span_\R \{ u,[\omega_0]\}\cap H^2 (M;\Z) = \Span_\Z \{ u\}.
\]
Assume by contradiction that this is false. Then $[\omega_0]=\kappa_1 u + \kappa_2 x$ for some $\kappa_1,\kappa_2\in\R$ and $x\in u^\bot \cap H^2 (M;\Z)$. Since $H^2 (M;\Z)$ together with the intersection pairing is an even unimodular quadratic lattice (see part (B2) of \ref{_unimodular-lattices-divisor-of-any-integral-vector-is-one_Proposition_}),
the rank of $\big(\Span_\R \{ u,x\}\big)^\bot \cap H^2 (M;\Z)$ equals $rk\, H^2 (M;\Z) - 2 = 20$. On the other hand, $\big(\Span_\R \{ u,x\}\big)^\bot \cap H^2 (M;\Z)$ is contained in $[\omega_0]^\bot\cap H^2 (M;\Z)$, whose rank is $7$,
yielding a contradiction. Thus, $\Span_\R \{ u,[\omega_0]\}\cap H^2 (M;\Z) = \Span_\Z \{ u\}$ and the claim holds.

Since, as we have shown, $[\omega_0]$ is orthoisotropically irrational,
\ref{_dense-orbit-of-sympl-form-K3-case_Theorem_} implies that the $\Diff^+ (M)$-orbit of $\omega_0$ is dense in $\Subot$ for each isotropic $u\in H^2 (M;\Z)$.

Now for each such $u$ denote by $\cLu$ the (possibly empty) set of $\omega\in\SymplKone$ such that $(M,\omega)$ admits a Maslov-zero Lagrangian torus in the homology class Poincar\'e-dual to $u$.
Clearly, $\cLu\subset \Subot$. Moreover, it follows from \ref{_Deforming-zero-homologous-Lagr-submfds_Proposition_} that $\cLu$ is open in $\Subot$.

The claim of the theorem is equivalent to showing that $\cLu$ is empty if $u$ is zero or non-primitive.

Assume, by contradiction, that $\cLu\neq\emptyset$ for such a $u$.
Hence, by the density property above, $\Diff^+ (M)$-orbit of $\omega_0$ intersects $\cLu$ -- that is, there exists $\omega$ in the $\Diff^+ (M)$-orbit of $\omega_0$ such that $(M,\omega)$ admits a Maslov-zero Lagrangian torus in the zero or a non-primitive homology class.
Then, by the $\Diff^+ (M)$-invariance, so does $(M,\omega_0)$, which yields a contradiction with the choice of $\omega_0$.

This finishes the proof of the theorem.
\endproof


\hfill


\noindent
{\bf Acknowledgments:} We thank N.Sheridan for pointing out a mistake in an earlier version of this paper.
 We also thank the anonymous referees for remarks and corrections.


\hfill


\bigskip

{\small
\noindent {\sc Michael Entov\\
Department of Mathematics\\
Technion - Israel Institute of Technology\\
Haifa 32000, Israel}\\
{\tt  entov@technion.ac.il}
}
\\

{\small
\noindent {\sc Misha Verbitsky\\
            {\sc Instituto Nacional de Matem\'atica Pura e
              Aplicada (IMPA) \\ Estrada Dona Castorina, 110\\
Jardim Bot\^anico, CEP 22460-320\\
Rio de Janeiro, RJ - Brasil}

\smallskip
\noindent
and

\smallskip
\noindent
Laboratory of Algebraic Geometry\\
National Research University HSE, Faculty of Mathematics\\
6 Usacheva Str., Moscow, Russia}\\
{\tt  verbit@impa.br}
}

\end{document}